\documentclass[11pt,twoside,notitlepage]{book}%
\usepackage{graphicx}
\usepackage{amsmath}
\usepackage[twoside]{geometry}
\usepackage{amsfonts}
\usepackage{amssymb}%
\newtheorem{theorem}{Theorem}
\newtheorem{acknowledgement}[theorem]{Acknowledgement}

\newtheorem{condition}[theorem]{Condition}

\newtheorem{corollary}[theorem]{Corollary}

\newtheorem{definition}[theorem]{Definition}
\newtheorem{example}[theorem]{Example}

\newtheorem{lemma}[theorem]{Lemma}

\newtheorem{proposition}[theorem]{Proposition}
\newtheorem{remark}[theorem]{Remark}

\newenvironment{proof}[1][Proof]{\textbf{#1.} }{\ \bf{q.e.d.}}
\begin{document}

\frontmatter
\title{New Li--Yau--Hamilton Inequalities for the Ricci Flow via the Space-time Approach}
\author{Bennett Chow\\University of California, San Diego
\and Dan Knopf\\University of Wisconsin}
\date{March 5, 2000 (revised April 8, 2002)}
\maketitle

\mainmatter

\chapter{Introduction}

In \cite{HamHarnack}, Hamilton determined a sharp differential Harnack
inequality of Li--Yau type for complete solutions of the Ricci flow with
non-negative curvature operator. This Li--Yau--Hamilton inequality
(abbreviated as \emph{LYH inequality }below) is of critical importance to the
understanding of singularities of the Ricci flow, as is evident from its
numerous applications in \cite{HamSurf}, \cite{HamEternal}, \cite{HamForm},
and \cite{HamPIC}. Moreover, it has been informally claimed by Hamilton that
the discovery of a LYH inequality in dimension $3$ valid without any
hypothesis on curvature is the main unresolved step in his program of
approaching Thurston's Geometrization Conjecture by applying the Ricci flow to
closed $3$-manifolds. See \cite{HamForm} for some of the reasons why such an
inequality is believed to hold. (One may also consult the survey paper
\cite{CaoChow}.) Based on unpublished research of Hamilton and Hamilton--Yau,
the search for such a LYH inequality appears to be an extremely complex and
delicate problem. Roughly speaking, their approach is to start with the
$3$-dimensional LYH inequality for solutions with nonnegative sectional
curvature and try to perturb that estimate so that it holds for solutions with
arbitrary initial data. Because of an estimate of Hamilton \cite{HamForm} and
Ivey \cite{Ivey} which shows that the curvature operator of $3$-dimensional
solutions tends in a sense to become nonnegative, there is hope that such a
procedure will work. Some unpublished work of Hamilton and Yau appears close
to establishing a general LYH inequality in dimension $3.$ However, so far no
such inequality is known.

Due to the perturbational nature of the existing approaches, it is also of
interest to understand how general a LYH inequality one can prove under the
original hypothesis of nonnegative curvature operator. In this direction,
Hamilton and one of the authors \cite{CH} obtained a linear trace LYH
inequality for a system consisting of a solution of the Lichnerowicz-Laplacian
heat equation for symmetric $2$-tensors coupled to a solution of the Ricci
flow. Since the pair of the Ricci and metric tensors of a solution to the
Ricci flow forms such a system, their linear trace inequality generalizes the
traced case of Hamilton's tensor (matrix) inequality. In \cite{HamSurf}
Hamilton had already observed the formal similarity between his proof of the
$2$-dimensional trace LYH inequality for the Ricci flow and Li and Yau's proof
\cite{LY} of their Harnack inequality for the heat equation on Riemannian
manifolds. In a sense, the linear trace inequality generalizes this link to
higher dimensions. In dimension $2$, meanwhile, the link was made stronger and
more evident by the discovery \cite{Chow} of a $1$-parameter family of
inequalities interpolating between the Li--Yau and linear trace estimates.

In another direction, one recalls that Hamilton's matrix inequality is
equivalent to the positivity of a certain quadratic form. Hamilton observed
\cite{HamForm} that the evolution equation satisfied by that quadratic
suggests that his LYH inequality may be some sort of extension of nonnegative
curvature operator. This was shown to be true by S.-C. Chu and one of the
authors in \cite{ChowChu1}. They introduced a degenerate metric and a certain
compatible connection on space and time that extends the Levi-Civita
connection of a solution of the Ricci flow. They proved that Hamilton's
quadratic is in fact the curvature of that connection. Because the space-time
metric and connection satisfy the Ricci flow for degenerate metrics, one can
then apply the methods of \cite{Ham3man} to show that the quadratic satisfies
a nice evolution equation. This fact is the starting point for the present paper.

In this paper, we prove a new differential Harnack inequality of
Li--Yau--Hamilton type for the Ricci flow by generalizing the construction in
\cite{ChowChu1}. Our inequality applies to solutions of the Ricci flow coupled
to a $1$-form and a $2$-form solving Hodge-Laplacian heat-type equations. In
this sense, one may regard it as a linear-type matrix LYH inequality. In its
general form, it looks quite different from Hamilton's matrix inequality ---
except in the K\"{a}hler case, where as a special case, one may take the
$1$-form to be the exterior derivative of the scalar curvature and the
$2$-form to be the Ricci form, thereby obtaining an inequality slightly weaker
but qualitatively similar to Hamilton's. (Note that Cao \cite{Cao} has
extended Hamilton's LYH inequality in the K\"{a}hler case to solutions with
nonnegative bisectional curvature.) We state the general form of our result
(Theorem \ref{main}) as our\medskip

\noindent\textbf{Main Theorem: }\emph{Let }$\left(  \mathcal{M},g\left(
t\right)  \right)  $\emph{ be a solution of the Ricci flow on a closed
manifold and a time interval }$[0,\Omega)$\emph{. Let }$A_{0}$\emph{ be a }%
$2$\emph{-form which is closed at }$t=0$\emph{, and let }$E_{0}$\emph{ be a
}$1$\emph{-form which is closed at }$t=0$\emph{. Then there is a solution
}$A\left(  t\right)  $\emph{ of}%
\[
\frac{\partial}{\partial t}A=\Delta_{d}A,\quad\quad\quad\quad A\left(
0\right)  =A_{0}%
\]
\emph{and a solution }$E\left(  t\right)  $\emph{ of}%
\[
\frac{\partial}{\partial t}E=\Delta_{d}E-d\left\vert A\right\vert _{g}%
^{2},\quad\quad\quad\quad E\left(  0\right)  =E_{0}%
\]
\emph{which exist for }$t\in\lbrack0,\Omega)$\emph{, where} $-\Delta
_{d}\doteqdot d\delta+\delta d$\emph{ is the Hodge--de Rham Laplacian. Suppose
that the quadratic}%
\begin{align*}
&  \Psi\left(  A,E,U,W\right) \\
&  \doteqdot\operatorname*{Rm}\left(  U,U\right)  -2\left\langle \nabla
_{W}A,U\right\rangle +\left\vert A\left(  W\right)  \right\vert ^{2}%
-\left\langle \nabla_{W}E,W\right\rangle \\
&  =R_{ijk\ell}U^{ij}U^{\ell k}+2W^{j}\nabla_{j}A_{k\ell}U^{\ell k}+\left(
g^{pq}A_{jp}A_{\ell q}-\nabla_{j}E_{\ell}\right)  W^{j}W^{\ell}%
\end{align*}
\emph{is non-negative at }$t=0$\emph{ for any }$2$\emph{-form }$U$\emph{ and
}$1$\emph{-form }$W$\emph{. Then the matrix inequality }$\Psi\left(
A,E,U,W\right)  \geq0$\emph{ persists for all }$t\in\lbrack0,\Omega
)$\emph{.\medskip}

The above linear-type inequality is a special case of Theorem \ref{full}
obtained by taking a limit which actually scales away part of the main highest
order terms in the more general LYH matrix inequality established in Theorem
\ref{full}.\medskip

\noindent\textbf{Corollary A: }\emph{Under the hypotheses above, the trace
inequality}%
\[
0\leq\psi\left(  A,E,W\right)  \doteqdot\operatorname*{Rc}\left(  W,W\right)
-2\left(  \delta A\right)  \left(  W\right)  +\left\vert A\right\vert
^{2}+\delta E
\]
\emph{persists for all }$t\in\lbrack0,\Omega)$\emph{.}

\medskip

\noindent\textbf{Corollary B: }\emph{Let }$\left(  \mathcal{M},g\left(
t\right)  \right)  $ \emph{be a K\"{a}hler solution of the Ricci flow with
non-negative curvature operator on a closed manifold $\mathcal{M}$. Then for
any }$2$\emph{-form }$U$\emph{, }$1$\emph{-form }$W$\emph{, and all }%
$t>0$\emph{ such that the solution exists, one has the matrix estimate}%
\begin{align*}
0  &  \leq\operatorname*{Rm}\left(  U,U\right)  -2\left\langle \nabla_{W}%
\rho,U\right\rangle +\frac{1}{4t^{2}}\left\vert W\right\vert ^{2}+\frac{1}%
{t}\operatorname*{Rc}\left(  W,W\right) \\
&  +\operatorname*{Rc}{}^{2}\left(  W,W\right)  +\frac{1}{2}\left(
\nabla\nabla R\right)  \left(  W,W\right)  ,
\end{align*}
\emph{where }$\rho$\emph{ is the Ricci form. By setting }$U=X\wedge W$\emph{
and tracing over }$W$\emph{, this implies the trace estimate}%
\[
0\leq\frac{\partial}{\partial t}R+\frac{2R}{t}+\frac{n}{2t^{2}}+2\left\langle
\nabla R,X\right\rangle +2\operatorname*{Rc}\left(  X,X\right)
\]
\emph{for any }$1$\emph{-form }$X$.\medskip

Although this LYH inequality is weaker than the trace inequality special case
of the matrix inequality in \cite{Cao}, it qualitatively similar, and it
arises from a much more general inequality. \medskip

\textbf{Corollary C: }\emph{Let }$\left(  \mathcal{M}^{2},g\left(  t\right)
\right)  $\emph{ be a solution of the Ricci flow on a closed surface. If
}$\left(  \phi,f\right)  $\emph{ is a pair solving the system}%
\begin{align*}
\frac{\partial}{\partial t}\phi &  =\Delta\phi+R\phi\\
\frac{\partial}{\partial t}f  &  =\Delta f+\phi^{2},
\end{align*}
\emph{then the trace inequality}%
\[
0\leq R\left|  X\right|  ^{2}+2\left\langle \nabla\phi,X\right\rangle
+\frac{\partial}{\partial t}f
\]
\emph{is preserved.}

\bigskip

This paper is structured as follows:

\begin{itemize}
\item In \S \ref{Glossary}, we extend the methods of \cite{ChowChu1} to the
case of the Ricci flow with a cosmological term $\mu$. Only by doing so for
$\mu\equiv1/2$ are we able to display Hamilton's differential Harnack
quadratic of Li-Yau type as exactly equal to the curvature of a space-time
connection, and thus to provide the reader with a precise glossary between the
space-time approach and the computations in \cite{HamHarnack}. A similar but
less precise correspondence was earlier established in \cite{ChowChu1}.

\item In \S \ref{Generalize}, we study all symmetric space-time connections
that are compatible with the degenerate space-time metric and evolve via the
Ricci flow for degenerate metrics. Because these connections are not unique,
their curvature tensors yield new Li--Yau--Hamilton inequalities, which
include Hamilton's matrix inequality as a special case. We then employ scaling
arguments to derive a non-negative symmetric bilinear form on space-time,
which is equivalent to the quadratic $\Psi$ described in classical language in
the Main Theorem, and whose traced form yields Corollary A.

\item In \S \ref{Apply}, we develop some examples in order to compare a few
special cases of the new Li--Yau--Hamilton inequality with known results. In
particular, we derive Corollary B (Proposition \ref{full-kaehler}) and
Corollary C (a result of tracing the matrix inequality in Proposition
\ref{surface-system}).
\end{itemize}

\begin{acknowledgement}
The authors wish to thank Grigori Perelman for calling their attention to an
error in the original preprint of this paper.
\end{acknowledgement}

\chapter[The Ricci flow rescaled]{The Ricci flow rescaled by a cosmological
term\label{Glossary}}

\section[Self-similar solutions]{Self-similar solutions of the Ricci flow}

In this section, we recall the equations for self-similar solutions to the
Ricci flow (called \emph{Ricci solitons }by Hamilton) in order to motivate the
introduction of the Ricci flow with a cosmological term. The basic reference
is \S 3 of \cite{HamHarnack}.

\begin{definition}
A solution $\left(  \mathcal{M}^{n},g\left(  t\right)  \right)  $ of the Ricci
flow%
\[
\frac{\partial}{\partial t}g=-2\operatorname*{Rc}%
\]
on a time interval $\mathcal{I}$ containing $0$ is called a \emph{homothetic
Ricci soliton }if%
\begin{equation}
g\left(  \cdot,t\right)  =a\left(  t\right)  \left(  \phi_{t}^{\ast}\hat
{g}\right)  \left(  \cdot\right)  \label{soliton-equation1}%
\end{equation}
for some fixed metric $\hat{g}$ on $\mathcal{M}$, some function $a$ of time
satisfying $a\left(  0\right)  =1$, and some $1$-parameter family of
diffeomorphisms $\left\{  \phi_{t}:t\in\mathcal{I}\right\}  $ generated by
vector fields $-V\left(  t\right)  $ with the property that their dual
$1$-forms are closed:%
\begin{equation}
\nabla_{i}V_{j}=\nabla_{j}V_{i}. \label{closed-1-forms}%
\end{equation}
In this case, we say that $\left(  \mathcal{M}^{n},g\left(  t\right)  \right)
$ \emph{flows along} $V$.
\end{definition}

It is not obvious that the representation (\ref{soliton-equation1}) is unique,
because the family $\left\{  \phi_{t}\right\}  $ may contain homotheties. We
put equation (\ref{soliton-equation1}) in a canonical form as follows:

\begin{lemma}
Suppose $g$ is a homothetic Ricci soliton having the form
(\ref{soliton-equation1}). Let $\hat{a}\doteqdot\dot{a}\left(  0\right)  $,
and let $\left\{  \psi_{t}:t\in\mathcal{I}\right\}  $ be the $1$-parameter
family of diffeomorphisms generated by the vector fields%
\[
-\frac{1}{1+\hat{a}t}V\left(  0\right)  ,
\]
with $\psi_{0}=\operatorname*{id}_{\mathcal{M}}$. Then%
\[
g\left(  \cdot,t\right)  =\left(  1+\hat{a}t\right)  \left(  \psi_{t}^{\ast
}\hat{g}\right)  \left(  \cdot\right)  .
\]

\end{lemma}

\begin{proof}
Let $G\left(  t\right)  $ be a smooth $1$-parameter family of metrics, and let
$\left\{  \theta_{t}\right\}  $ be a family of diffeomorphisms generated by
vector fields $-W\left(  t\right)  $. Then we have%
\begin{align}
\frac{\partial}{\partial t}\left(  \theta_{t}^{\ast}G\left(  t\right)
\right)   &  =\left.  \frac{\partial}{\partial s}\right\vert _{s=0}\left(
\theta_{t+s}^{\ast}G\left(  t+s\right)  \right) \nonumber\\
&  =\left.  \theta_{t}^{\ast}\left(  \frac{\partial}{\partial t}G\left(
t\right)  \right)  +\frac{\partial}{\partial s}\right\vert _{s=0}\left[
\left(  \theta_{t}^{-1}\circ\theta_{t+s}\right)  ^{\ast}\left(  \theta
_{t}^{\ast}G\left(  t\right)  \right)  \right] \nonumber\\
&  =\theta_{t}^{\ast}\left(  \frac{\partial}{\partial t}G\left(  t\right)
\right)  -\mathcal{L}_{\left(  \theta_{t}^{-1}\right)  _{\ast}W\left(
t\right)  }\left(  \theta_{t}^{\ast}G\left(  t\right)  \right)  .\nonumber\\
&  =\theta_{t}^{\ast}\left(  \frac{\partial}{\partial t}G\left(  t\right)
-\mathcal{L}_{W\left(  t\right)  }G\left(  t\right)  \right)  .
\label{DifferentiatePullback}%
\end{align}
Applying (\ref{DifferentiatePullback}) with $G\left(  t\right)  =a\left(
t\right)  \hat{g}$ and $\theta_{t}=\phi_{t}$, we get%
\begin{equation}
\operatorname{Rc}\left(  g\right)  =-\frac{1}{2}\frac{\partial}{\partial
t}g\left(  t\right)  =-\frac{1}{2}\frac{\dot{a}}{a}g+\frac{a}{2}\phi_{t}%
^{\ast}\left(  \mathcal{L}_{V\left(  t\right)  }\hat{g}\right)  .
\label{soliton-equation2}%
\end{equation}
Now define%
\[
\tilde{g}\left(  t\right)  \doteqdot\left(  1+\hat{a}t\right)  \left(
\psi_{t}^{\ast}\hat{g}\right)  .
\]
Applying (\ref{DifferentiatePullback}) with $G\left(  t\right)  =\left(
1+\hat{a}t\right)  \hat{g}$ and $\theta_{t}=\psi_{t}$, we obtain%
\[
\frac{\partial}{\partial t}\tilde{g}\left(  t\right)  =\psi_{t}^{\ast}\left(
\hat{a}\hat{g}-\mathcal{L}_{\frac{1}{1+\hat{a}t}V\left(  0\right)  }\left(
\left(  1+\hat{a}t\right)  \hat{g}\right)  \right)  =\psi_{t}^{\ast}\left(
\hat{a}\hat{g}-\mathcal{L}_{V\left(  0\right)  }\hat{g}\right)  .
\]
But evaluating equation (\ref{soliton-equation2}) at $t=0$ shows that%
\[
\operatorname{Rc}\left(  \hat{g}\right)  =\left.  \operatorname{Rc}\left(
g\right)  \right\vert _{t=0}=-\frac{1}{2}\hat{a}\hat{g}+\frac{1}{2}%
\mathcal{L}_{V\left(  0\right)  }\hat{g}.
\]
Hence%
\[
\frac{\partial}{\partial t}\tilde{g}\left(  t\right)  =\psi_{t}^{\ast}\left(
-2\operatorname{Rc}\left(  \hat{g}\right)  \right)  =-2\operatorname{Rc}%
\left(  \tilde{g}\left(  t\right)  \right)  .
\]
So $\tilde{g}\left(  t\right)  $ is a solution of the Ricci flow with
$\tilde{g}\left(  0\right)  =g\left(  0\right)  $. Since solutions of the
Ricci flow are unique, it follows that $\tilde{g}\left(  t\right)  =g\left(
t\right)  $ for as long as both solutions exist.
\end{proof}

The gist of the lemma is that in (\ref{soliton-equation1}) and
(\ref{soliton-equation2}) we may assume $\dot{a}\left(  t\right)  \equiv
\hat{a}$ is independent of $t$, so that $a\left(  t\right)  =1+\hat{a}t$ and
$a\left(  t\right)  V\left(  t\right)  =V\left(  0\right)  .$

From the point of view of motivating the differential Harnack inequality of
Li-Yau type, Hamilton considered the extreme case for the function $a\left(
t\right)  $ in the definition of Ricci soliton. In particular, he was
interested in the case where the initial metric $g\left(  0\right)  $ is
singular (such as the metric of a cone) and the metric $g\left(  t\right)  $
expands as $t$ increases. Formally, this corresponds to letting $\hat{a}%
=\dot{a}\left(  0\right)  $ tend to infinity,%
\begin{equation}
\lim_{\hat{a}\rightarrow\infty}\frac{\dot{a}}{a}=\lim_{\hat{a}\rightarrow
\infty}\frac{\hat{a}}{1+\hat{a}t}=\frac{1}{t}, \label{cone}%
\end{equation}
and motivates the following definition:

\begin{definition}
\label{ExpRcSoliton}A solution $\left(  \mathcal{M}^{n},g\left(  t\right)
\right)  $ of the Ricci flow%
\[
\frac{\partial}{\partial t}g=-2\operatorname*{Rc}%
\]
on a time interval $\left(  0,\Omega\right)  $ containing $1$ is called an
\emph{expanding Ricci soliton flowing along }$V$\emph{ }if%
\begin{equation}
g\left(  \cdot,t\right)  =t\left(  \theta_{t}^{\ast}\hat{g}\right)  \left(
\cdot\right)  \label{exp-ric-sol-eqn1}%
\end{equation}
for some fixed metric $\hat{g}$ on $\mathcal{M}$ and some $1$-parameter family
of diffeomorphisms $\left\{  \theta_{t}:t\in\left(  0,\Omega\right)  \right\}
$ such that $\theta_{1}=\operatorname*{id}_{\mathcal{M}}$ and the dual
$1$-forms of the vector fields $-V\left(  t\right)  $ which generate
$\theta_{t}$ are closed.
\end{definition}

Differentiating (\ref{exp-ric-sol-eqn1}) leads to the equation%
\begin{equation}
R_{ij}=\nabla_{i}V_{j}-\frac{1}{2t}g_{ij}, \label{exp-ric-sol-eqn2}%
\end{equation}
where $V\left(  t\right)  =\frac{1}{t}V\left(  1\right)  $. Notice that
(\ref{exp-ric-sol-eqn2}) can be obtained formally by passing to the limit
(\ref{cone}) in equation (\ref{soliton-equation2}).

\bigskip

Now observe that $t^{-1}g\left(  t\right)  =\theta_{t}^{\ast}\left(  g\left(
1\right)  \right)  $ evolves by diffeomorphisms. This motivates us to make the
following transformation for any solution $g\left(  t\right)  $ to the Ricci
flow:%
\[
\bar{g}\left(  t\right)  \doteqdot\frac{1}{t}g\left(  t\right)  .
\]
To get a nice equation for $\bar{g}$, we change the time variable by $\bar
{t}\doteqdot\ln t.$ Then $\bar{g}\left(  \bar{t}\right)  $ is a solution to
the equation%
\[
\frac{\partial}{\partial\bar{t}}\bar{g}_{ij}=-2\left(  \bar{R}_{ij}+\frac
{1}{2}\bar{g}_{ij}\right)
\]
on the time interval $\left(  -\infty,\ln\Omega\right)  $ containing $\bar
{t}=0$. We call this equation the \emph{Ricci flow with cosmological constant}
$1/2$. More generally, there is the following:

\begin{definition}
We say that $\left(  \mathcal{M}^{n},g\left(  \bar{t}\right)  \right)  $ is a
solution of the \emph{Ricci flow with cosmological term} $\mu\left(  \bar
{t}\right)  \in\mathbb{R}$ on a time interval $\mathcal{I}$ if%
\begin{equation}
\frac{\partial}{\partial\bar{t}}\bar{g}_{ij}=-2\left(  \bar{R}_{ij}+\mu\bar
{g}_{ij}\right)  . \label{cosmo}%
\end{equation}

\end{definition}

A solution of the Ricci flow with cosmological constant $\mu=1/2$ is an
expanding Ricci soliton%
\[
-\frac{1}{2}\frac{\partial}{\partial t}g_{ij}=R_{ij}=\nabla_{i}V_{j}-\frac
{1}{2t}g_{ij}%
\]
if and only if $\bar{g}\left(  \bar{t}\right)  =e^{-\bar{t}}g\left(
e^{\bar{t}}\right)  $ satisfies%
\begin{equation}
-\frac{1}{2}\frac{\partial}{\partial\bar{t}}\bar{g}_{ij}=\bar{R}_{ij}+\frac
{1}{2}\bar{g}_{ij}=\bar{\nabla}_{i}\bar{V}_{j}, \label{cosmo-soliton-equation}%
\end{equation}
where $\bar{V}_{j}\left(  \bar{t}\right)  \doteqdot V_{j}\left(  t\right)  $,
hence if and only if $\bar{g}\left(  \bar{t}\right)  $ is a steady Ricci
soliton. Note that $\bar{V}^{k}\left(  \bar{t}\right)  =\bar{g}^{jk}\left(
\bar{t}\right)  \bar{V}_{j}\left(  \bar{t}\right)  =tV^{k}\left(  t\right)
=V^{k}\left(  1\right)  $ is independent of $\bar{t}$, so that
\begin{equation}
\frac{\partial}{\partial\bar{t}}\bar{V}^{k}=0.
\label{cosmo-soliton-vector-time-indep}%
\end{equation}
Taking the divergence of (\ref{cosmo-soliton-equation}), using
(\ref{closed-1-forms}), and commuting derivatives, we compute%
\[
\frac{1}{2}\bar{\nabla}_{j}\bar{R}=\bar{\nabla}_{i}\bar{\nabla}^{i}\bar{V}%
_{j}=\bar{\nabla}_{j}\bar{\nabla}_{i}\bar{V}^{i}+\bar{R}_{jk}\bar{V}^{k}%
=\bar{\nabla}_{j}\left(  \bar{R}+\frac{n}{2}\right)  +\bar{R}_{jk}\bar{V}^{k}%
\]
and thus obtain the following useful identities, valid when $\mu\equiv\frac
{1}{2}$:%
\begin{equation}
\frac{1}{2}\bar{\nabla}\bar{R}=\bar{\Delta}\bar{V}=-\overline
{\operatorname*{Rc}}\left(  \bar{V}\right)  . \label{cosmo-soliton-div-eqn}%
\end{equation}

\section[The space-time connection]{The space-time connection for the Ricci
flow rescaled by a cosmological term}

In this section, we show that the definition of the space-time connection for
the Ricci flow in \cite{ChowChu1} may be extended to the case where there is a
cosmological term $\mu\left(  \bar{t}\right)  $ in the flow equation.
Motivated by the discussion in the previous section, we are mainly interested
in the case $\mu\equiv\frac{1}{2}$. Since the relevant computations are
modifications of those in \cite{ChowChu1}, we shall omit many details of the proofs.

Let $\mathcal{\ }\widetilde{\mathcal{M}}=\mathcal{M}^{n}\times\mathcal{I}$ and
denote the time coordinate by $x^{0}\doteqdot\bar{t}$. Recall that the
degenerate space-time metric on $T^{\ast}\widetilde{\mathcal{M}}$ is defined
by%
\begin{equation}
\tilde{g}^{ij}\doteqdot\left\{
\begin{tabular}
[c]{ll}%
$\bar{g}^{ij}$ & if $i,j\geq1$\\
$0$ & if $i=0$ or $j=0$%
\end{tabular}
\right.  . \label{define-gtilde}%
\end{equation}
Modifying the definition in \cite{ChowChu1}, we define a symmetric space-time
connection $\tilde{\nabla}$ by specifying its Christoffel symbols to be%
\begin{align}
\tilde{\Gamma}_{ij}^{k}  &  =\bar{\Gamma}_{ij}^{k}\tag{C1}\label{C1}\\
\tilde{\Gamma}_{i0}^{k}  &  =-\left(  \bar{R}_{i}^{k}+\mu\delta_{i}^{k}\right)
\tag{C2}\label{C2}\\
\tilde{\Gamma}_{00}^{k}  &  =-\frac{1}{2}\bar{\nabla}^{k}\bar{R}%
\tag{C3}\label{C3}\\
\tilde{\Gamma}_{00}^{0}  &  =-\mu\tag{C4}\label{C4}\\
\tilde{\Gamma}_{ij}^{0}  &  =\tilde{\Gamma}_{i0}^{0}=0, \tag{C5}\label{C5}%
\end{align}
where $i,j,k\geq1$.

\begin{lemma}
\label{st-metric-parallel}The connection $\tilde{\nabla}$ is compatible with
the degenerate metric $\tilde{g}$: for all $i,j,k\geq0$,
\[
\tilde{\nabla}_{i}\tilde{g}^{jk}=0.
\]

\end{lemma}

\begin{proof}
This is a straightforward computation using the identity
\[
\tilde{\nabla}_{i}\tilde{g}^{jk}=\partial_{i}\tilde{g}^{jk}+\tilde{\Gamma
}_{ip}^{j}\tilde{g}^{pk}+\tilde{\Gamma}_{ip}^{k}\tilde{g}^{jp}%
\]
with formulas (\ref{C1}), (\ref{C2}), and (\ref{C5}).
\end{proof}

Given a time-dependent vector field $\bar{W}\left(  \bar{t}\right)  $ on
$\mathcal{M}$, we associate to it the space-time vector field%
\[
\tilde{W}\left(  \bar{t}\right)  \doteqdot\frac{\partial}{\partial\bar{t}%
}+\bar{W}\left(  \bar{t}\right)  .
\]
In local coordinates, $\tilde{W}^{0}=1$ and $\tilde{W}^{j}=\bar{W}^{j}$ if
$j\geq1$.

\begin{lemma}
The formulas for the covariant derivative of the vector field $\tilde{W}$ are%
\begin{align}
\tilde{\nabla}_{i}\tilde{W}^{j}  &  =\bar{\nabla}_{i}\bar{W}^{j}-\left(
\bar{R}_{i}^{j}+\mu\delta_{i}^{j}\right) \tag{CW1}\label{P1}\\
\tilde{\nabla}_{0}\tilde{W}^{j}  &  =\frac{\partial}{\partial\bar{t}}\bar
{W}^{j}-\left(  \bar{R}_{k}^{j}+\mu\delta_{k}^{j}\right)  \bar{W}^{k}-\frac
{1}{2}\bar{\nabla}^{j}\bar{R}\tag{CW2}\label{P2}\\
\tilde{\nabla}_{0}\tilde{W}^{0}  &  =-\mu\tag{CW3}\label{P3}\\
\tilde{\nabla}_{i}\tilde{W}^{0}  &  =0 \tag{CW4}\label{P4}%
\end{align}
for all $i,j,k\geq1$.
\end{lemma}

\begin{proof}
This follows from $\tilde{\nabla}_{i}\tilde{W}^{j}=\partial_{i}\tilde{W}%
^{j}+\sum_{p=1}^{n}\tilde{\Gamma}_{ip}^{j}\tilde{W}^{p}+\tilde{\Gamma}%
_{i0}^{j}\tilde{W}^{0}$ and all of the formulas (\ref{C1})--(\ref{C5}).
\end{proof}

We can now make the important observation that \emph{the space-time of a
steady soliton flowing along }$\bar{V}$\emph{ has a geometric product
structure. }Recall that a parallel vector field on a Riemannian manifold
$\mathcal{N}$ gives a local splitting of $\mathcal{N}$ as the product of an
open interval with an $\left(  n-1\right)  $-dimensional manifold
$\mathcal{P}$. Hence the observation follows from:

\begin{proposition}
If $\bar{g}\left(  \bar{t}\right)  $ is a steady soliton of the Ricci flow
with cosmological constant $\mu=\frac{1}{2}$ flowing along the vector fields
$\bar{V}\left(  \bar{t}\right)  $, then
\[
\tilde{\nabla}_{i}\left(  e^{\frac{1}{2}\bar{t}}\tilde{V}\right)  ^{j}\equiv0
\]
for all $i,j\geq0$. That is, the space-time vector field $e^{\frac{1}{2}%
\bar{t}}\tilde{V}$ is parallel.
\end{proposition}

\begin{proof}
If $i=j=0$, the formula follows from (\ref{P3}). For the case $i=0,\;j\geq1$,
we apply equations (\ref{P2}), (\ref{cosmo-soliton-vector-time-indep}), and
(\ref{cosmo-soliton-div-eqn}). If $i\geq1$ and $j=0$, the formula follows from
(\ref{P4}). And for the case $i\geq1,\;j\geq1$, we apply (\ref{P1}) and
(\ref{cosmo-soliton-equation}).
\end{proof}

\subsection{The Riemann curvature tensor}

Denote the Riemann curvature tensor of the space-time connection
$\tilde{\nabla}$ by
\begin{equation}
\tilde{R}\left(  \tilde{X},\tilde{Y}\right)  \tilde{Z}=\tilde{\nabla}%
_{\tilde{X}}\tilde{\nabla}_{\tilde{Y}}\tilde{Z}-\tilde{\nabla}_{\tilde{Y}%
}\tilde{\nabla}_{\tilde{X}}\tilde{Z}-\tilde{\nabla}_{\left[  \tilde{X}%
,\tilde{Y}\right]  }\tilde{Z}. \label{Rm-def}%
\end{equation}
Since $e^{\frac{1}{2}\bar{t}}\tilde{V}$ is parallel and $\tilde{R}\left(
\tilde{X},\tilde{Y}\right)  \tilde{Z}$ is a tensor, we immediately get:

\begin{corollary}
\label{Vparallel}If $\bar{g}\left(  \bar{t}\right)  $ is a steady soliton
flowing along $\bar{V}\left(  \bar{t}\right)  $ with $\mu=\frac{1}{2}$, then%
\begin{equation}
\tilde{R}\left(  \tilde{X},\tilde{Y}\right)  \tilde{V}=0
\label{cosmo-soliton-curvature-zero-dir}%
\end{equation}
for all $\tilde{X}$ and $\tilde{Y}.$
\end{corollary}

We shall see in the next section how this relates to the derivation of the
Li--Yau--Hamilton quadratic in \S 3 of \cite{HamPCO}.

\bigskip

The formulas for the space-time Riemann curvature tensor are as follows:

\begin{proposition}
\label{st-riem-formulas}If $i,j,k,\ell\geq1$ and $a,b,c\geq0$, then
$\widetilde{\operatorname*{Rm}}$ satisfies:%
\begin{align}
\tilde{R}_{ijk}^{\ell}  &  =\bar{R}_{ijk}^{\ell}\tag{R1}\label{R1}\\
\tilde{R}_{ij0}^{\ell}  &  =\bar{\nabla}_{j}\bar{R}_{i}^{\ell}-\bar{\nabla
}_{i}\bar{R}_{j}^{\ell}\tag{R2a}\label{R2a}\\
\tilde{R}_{0jk}^{\ell}  &  =\bar{\nabla}^{\ell}\bar{R}_{jk}-\bar{\nabla}%
_{k}\bar{R}_{j}^{\ell}\tag{R2b}\label{R2b}\\
\tilde{R}_{i00}^{\ell}  &  =\frac{\partial}{\partial\bar{t}}\bar{R}_{i}^{\ell
}-\frac{1}{2}\bar{\nabla}_{i}\bar{\nabla}^{\ell}\bar{R}-\bar{R}_{i}^{p}\bar
{R}_{p}^{\ell}-\mu\bar{R}_{i}^{\ell}+\frac{d\mu}{d\bar{t}}\delta_{i}^{\ell
}\tag{R3}\label{R3}\\
\tilde{R}_{abc}^{0}  &  =0. \tag{R4}\label{R4}%
\end{align}

\end{proposition}

\begin{remark}
The standard asymmetries satisfied by the curvature of any connection imply in
particular that%
\begin{align*}
\tilde{R}_{0jk}^{\ell}+\tilde{R}_{j0k}^{\ell}  &  =0\\
\tilde{R}_{i00}^{\ell}+\tilde{R}_{0i0}^{\ell}  &  =0.
\end{align*}
Because $\tilde{\nabla}$ is torsion-free, the first and second Bianchi
identities take the form:%
\begin{align}
\tilde{R}_{ijk}^{\ell}+\tilde{R}_{jki}^{\ell}+\tilde{R}_{kij}^{\ell}  &
=0\tag{B1}\label{B1}\\
\tilde{\nabla}_{m}\tilde{R}_{ijk}^{\ell}+\tilde{\nabla}_{i}\tilde{R}%
_{jmk}^{\ell}+\tilde{\nabla}_{j}\tilde{R}_{mik}^{\ell}  &  =0 \tag{B2}%
\label{B2}%
\end{align}
for all $i,j,k,\ell,m\geq0$.
\end{remark}

\begin{remark}
\label{R3-rewrite}Using the evolution equation%
\[
\frac{\partial}{\partial\bar{t}}\bar{R}_{i}^{\ell}=\bar{\Delta}\bar{R}%
_{i}^{\ell}+2\bar{R}_{ipq}^{\ell}\bar{R}^{pq}+2\mu\bar{R}_{i}^{\ell},
\]
we may rewrite (\ref{R3}) as%
\[
\tilde{R}_{i00}^{\ell}=\bar{\Delta}\bar{R}_{i}^{\ell}-\frac{1}{2}\bar{\nabla
}_{i}\bar{\nabla}^{\ell}\bar{R}+2\bar{R}_{ipq}^{\ell}\bar{R}^{pq}-\bar{R}%
_{i}^{p}\bar{R}_{p}^{\ell}+\mu\bar{R}_{i}^{\ell}+\frac{d\mu}{d\bar{t}}%
\delta_{i}^{\ell}.
\]

\end{remark}

The identities in Proposition (\ref{st-riem-formulas}) are proved in a manner
similar to Theorem 2.2 and 3.1 of \cite{ChowChu1}, which give the
corresponding equations for the case $\mu\equiv0$.

The components of the Ricci tensor are given by $\tilde{R}_{jk}=\Sigma
_{i=0}^{n}\tilde{R}_{ijk}^{i}=\Sigma_{i=1}^{n}\tilde{R}_{ijk}^{i}$. Hence
tracing (as in Corollary 2.4 of \cite{ChowChu1}) gives the following:

\begin{corollary}
The Ricci tensor satisfies the identities:%
\begin{align}
\tilde{R}_{ij}  &  =\bar{R}_{ij}\tag{Rc1}\label{Rc1}\\
\tilde{R}_{0j}  &  =\frac{1}{2}\bar{\nabla}_{j}\bar{R}\tag{Rc2}\label{Rc2}\\
\tilde{R}_{00}  &  =\frac{1}{2}\frac{\partial}{\partial\bar{t}}\bar{R}%
+n\frac{d\mu}{d\bar{t}}. \tag{Rc3}\label{Rc3}%
\end{align}

\end{corollary}

As in Lemma 3.3 of \cite{ChowChu1}, we notice that:

\begin{remark}
\label{cov-Ricci-sym}The covariant derivatives of the Ricci tensor obey the
symmetries%
\begin{align}
\tilde{\nabla}_{i}\tilde{R}_{j0}  &  =\tilde{\nabla}_{j}\tilde{R}%
_{i0}\tag{CRc1}\label{CRc1}\\
\tilde{\nabla}_{i}\tilde{R}_{00}  &  =\tilde{\nabla}_{0}\tilde{R}_{i0}
\tag{CRc2}\label{CRc2}%
\end{align}
for all $i,j\geq1.$
\end{remark}

\begin{proof}
Using (\ref{Rc2}) and (\ref{C2}), we find that%
\[
\tilde{\nabla}_{i}\tilde{R}_{j0}=\bar{\nabla}_{i}\bar{\nabla}_{j}\bar{R}%
+\bar{R}_{ij}^{2}+\mu\bar{R}_{ij}=\tilde{\nabla}_{j}\tilde{R}_{i0},
\]
which proves equation (\ref{CRc1}). Using (\ref{Rc2}) and (\ref{Rc3}), we get
\begin{align*}
\tilde{\nabla}_{i}\tilde{R}_{00}-\tilde{\nabla}_{0}\tilde{R}_{i0}  &
=\bar{\nabla}_{i}\left(  \frac{1}{2}\frac{\partial}{\partial\bar{t}}\bar
{R}\right)  -2\tilde{\Gamma}_{i0}^{p}\tilde{R}_{p0}\\
&  -\frac{\partial}{\partial\bar{t}}\left(  \frac{1}{2}\bar{\nabla}_{i}\bar
{R}\right)  +\tilde{\Gamma}_{0i}^{p}\tilde{R}_{p0}+\tilde{\Gamma}_{00}%
^{p}\tilde{R}_{ip}+\tilde{\Gamma}_{00}^{0}\tilde{R}_{i0}%
\end{align*}
where $p$ is summed from $1$ to $n$. Equation (\ref{CRc2}) then follows by
applying formulas (\ref{C2})--(\ref{C4}) and (\ref{Rc1})--(\ref{Rc3}).
\end{proof}

Generalizing the definition in \cite{ChowChu1}, we have the following:

\begin{definition}
\label{degenerateRF}A degenerate metric and compatible connection $\left(
\tilde{g},\tilde{\nabla}\right)  $ satisfy the Ricci flow with cosmological
term $\mu$ if for all $i,j,k\geq0$,%
\begin{equation}
\frac{\partial}{\partial\bar{t}}\tilde{\Gamma}_{ij}^{k}=\tilde{g}^{k\ell
}\left(  -\tilde{\nabla}_{i}\tilde{R}_{j\ell}-\tilde{\nabla}_{j}\tilde
{R}_{i\ell}+\tilde{\nabla}_{\ell}\tilde{R}_{ij}\right)  .
\label{cosmo-st-christoffel-time-deriv}%
\end{equation}

\end{definition}

\begin{proposition}
The pair $\left(  \tilde{g},\tilde{\nabla}\right)  $ satisfies the Ricci flow
with cosmological term $\mu$.
\end{proposition}

\begin{proof}
If $i,j,k\geq1$, the standard formula%
\[
\frac{\partial}{\partial\bar{t}}\bar{\Gamma}_{ij}^{k}=\frac{1}{2}\bar
{g}^{k\ell}\left[  \bar{\nabla}_{i}\left(  \frac{\partial}{\partial\bar{t}%
}\bar{g}_{\ell j}\right)  +\bar{\nabla}_{j}\left(  \frac{\partial}%
{\partial\bar{t}}\bar{g}_{i\ell}\right)  -\bar{\nabla}_{\ell}\left(
\frac{\partial}{\partial\bar{t}}\bar{g}_{ij}\right)  \right]
\]
shows that (\ref{cosmo-st-christoffel-time-deriv}) holds. If $k=0$, the result
is trivial by equations (\ref{C4}) and (\ref{C5}). If $i=0$ and $j,k\geq1$,
the observation $\tilde{\Gamma}_{0j}^{p}\tilde{R}_{p}^{k}=\tilde{\Gamma}%
_{op}^{k}\tilde{R}_{j}^{p}$ and identity (\ref{CRc1}) together imply that%
\[
-\tilde{\nabla}_{0}\tilde{R}_{j}^{k}-\tilde{\nabla}_{j}\tilde{R}_{0}%
^{k}+\tilde{\nabla}^{k}\tilde{R}_{0j}=-\tilde{\nabla}_{0}\tilde{R}_{j}%
^{k}=-\frac{\partial}{\partial\bar{t}}\bar{R}_{j}^{k}=\frac{\partial}%
{\partial\bar{t}}\bar{\Gamma}_{0j}^{k}.
\]
If $i=j=0$ and $k\geq1$, the observation $\tilde{\Gamma}_{00}^{p}\tilde{R}%
_{p}^{k}=\tilde{\Gamma}_{0p}^{k}\tilde{R}_{0}^{p}$ and identity (\ref{CRc2})
imply%
\[
-2\tilde{\nabla}_{0}\tilde{R}_{0}^{k}+\tilde{\nabla}^{k}\tilde{R}_{00}%
=-\tilde{\nabla}_{0}\tilde{R}_{0}^{k}=-\frac{\partial}{\partial\bar{t}}\left(
\frac{1}{2}\bar{\nabla}^{k}\bar{R}\right)  =\frac{\partial}{\partial\bar{t}%
}\tilde{\Gamma}_{00}^{k}.
\]

\end{proof}

\begin{lemma}
If $\mu$ is constant, the space-time curvature tensor satisfies the divergence
identity%
\begin{equation}
\tilde{g}^{pq}\tilde{\nabla}_{p}\tilde{R}_{qjk}^{\ell}=\tilde{R}_{0jk}^{\ell}
\label{st-curv-div}%
\end{equation}
between components of the $\left(  2,1\right)  $-tensor on the LHS and
components of the $\left(  3,1\right)  $-tensor on the RHS.
\end{lemma}

\begin{proof}
If $j\geq1$ and $k\geq1$, this is just the contracted second Bianchi identity%
\[
\tilde{g}^{pq}\tilde{\nabla}_{p}\tilde{R}_{qjk}^{\ell}=\bar{\nabla}^{p}\bar
{R}_{pjk}^{\ell}=\bar{\nabla}^{\ell}\bar{R}_{jk}-\bar{\nabla}_{k}\bar{R}%
_{j}^{\ell}.
\]
If $j\geq1$, $k=0$, and $\mu$ is constant, this follows from Remark
\ref{R3-rewrite} and the calculation%
\[
\tilde{g}^{pq}\tilde{\nabla}_{p}\tilde{R}_{qj0}^{\ell}=-\bar{\Delta}\bar
{R}_{j}^{\ell}+\frac{1}{2}\bar{\nabla}_{j}\bar{\nabla}^{\ell}\bar{R}-2\bar
{R}^{pq}\bar{R}_{jpq}^{\ell}+\bar{R}_{j}^{p}\bar{R}_{p}^{\ell}-\mu\bar{R}%
_{j}^{\ell}.
\]
If $j=0$ and $k\geq1$, one computes%
\[
\tilde{g}^{pq}\tilde{\nabla}_{p}\tilde{R}_{q0k}^{\ell}=\bar{\nabla}^{p}%
\bar{\nabla}_{k}\bar{R}_{p}^{\ell}-\bar{\nabla}^{p}\bar{\nabla}^{\ell}\bar
{R}_{pk}+\bar{R}^{pm}\bar{R}_{pmk}^{\ell}=2\bar{R}^{pm}\bar{R}_{pmk}^{\ell
}=0.
\]
Finally, if $j=k=0$, one obtains
\begin{align*}
\tilde{g}^{pq}\tilde{\nabla}_{p}\tilde{R}_{q00}^{\ell}  &  =\bar{\nabla}%
^{p}\bar{\Delta}\bar{R}_{p}^{\ell}-\frac{1}{2}\bar{\Delta}\bar{\nabla}^{\ell
}\bar{R}+2\bar{\nabla}^{p}\left(  \bar{R}_{pqm}^{\ell}\bar{R}^{qm}\right) \\
&  -\frac{1}{2}\bar{R}_{p}^{\ell}\bar{\nabla}^{p}\bar{R}-\bar{R}^{pq}%
\bar{\nabla}^{\ell}\bar{R}_{pq}\\
&  =0
\end{align*}
by a straightforward calculation. (See also Lemma 2.2 and the remark after it
in \cite{ChowChu2}).
\end{proof}

\begin{remark}
Tracing formula (\ref{st-curv-div}) and applying (\ref{B2}) yields%
\[
\tilde{R}_{0}^{\ell}=\tilde{g}^{jk}\tilde{g}^{pq}\tilde{\nabla}_{p}\tilde
{R}_{qjk}^{\ell}=\tilde{g}^{pq}\tilde{\nabla}_{p}\tilde{R}_{q}^{\ell}=\frac
{1}{2}\tilde{\nabla}^{\ell}\tilde{R}=\frac{1}{2}\bar{\nabla}^{\ell}\bar{R},
\]
in agreement with (\ref{Rc2}).
\end{remark}

The evolution equation for the space-time curvature tensor is given by:

\begin{proposition}
\label{st-riem-31-evol}If $\mu$ is constant, then%
\[
\tilde{\nabla}_{0}\tilde{R}_{ijk}^{\ell}=\tilde{\Delta}\tilde{R}_{ijk}^{\ell
}+2\left(  \tilde{B}_{ijk}^{\ell}-\tilde{B}_{jik}^{\ell}-\tilde{B}_{jki}%
^{\ell}+\tilde{B}_{ikj}^{\ell}\right)  +2\mu\tilde{R}_{ijk}^{\ell},
\]
where%
\[
\tilde{B}_{ijk}^{\ell}\doteqdot-\tilde{g}^{pq}\tilde{R}_{pij}^{m}\tilde
{R}_{kqm}^{\ell}.
\]

\end{proposition}

\begin{proof}
This formula may be proved along the lines of \cite{Ham3man}. Instead, we give
an alternate proof using the space-time Bianchi and divergence identities. We
note that taking the covariant derivative of identity (\ref{st-curv-div})
yields
\[
\tilde{\nabla}_{i}\tilde{R}_{0jk}^{\ell}=\tilde{\nabla}_{i}\left(  \tilde
{g}^{pq}\tilde{\nabla}_{p}\tilde{R}_{qjk}^{\ell}\right)  -\tilde{\Gamma}%
_{i0}^{m}\tilde{R}_{mjk}^{\ell}.
\]
So by using (\ref{B2}), substituting, and cancelling terms, we directly
obtain
\begin{align*}
\tilde{\nabla}_{0}\tilde{R}_{ijk}^{\ell}  &  =\tilde{\nabla}_{i}\tilde
{R}_{0jk}^{\ell}-\tilde{\nabla}_{j}\tilde{R}_{0ik}^{\ell}\\
&  =\tilde{g}^{pq}\left(  \tilde{\nabla}_{i}\tilde{\nabla}_{p}\tilde{R}%
_{qjk}^{\ell}-\tilde{\nabla}_{j}\tilde{\nabla}_{p}\tilde{R}_{qik}^{\ell
}\right)  -\tilde{\Gamma}_{i0}^{m}\tilde{R}_{mjk}^{\ell}+\tilde{\Gamma}%
_{j0}^{m}\tilde{R}_{mik}^{\ell}\\
&  =\tilde{\Delta}\tilde{R}_{ijk}^{\ell}+2\tilde{g}^{pq}\left(  \tilde
{R}_{ipk}^{m}\tilde{R}_{jqm}^{\ell}-\tilde{R}_{ipm}^{\ell}\tilde{R}_{jqk}%
^{m}\right) \\
&  -\tilde{g}^{pq}\tilde{R}_{ijp}^{m}\tilde{R}_{qmk}^{\ell}+2\mu\tilde
{R}_{ijk}^{\ell},
\end{align*}
where $\tilde{\Delta}\doteqdot\tilde{g}^{pq}\tilde{\nabla}_{p}\tilde{\nabla
}_{q}$ is the space-time Laplacian. Then using (\ref{B1}) and the identity
$\tilde{B}_{ijk}^{\ell}=-\tilde{g}^{pq}\tilde{R}_{pji}^{m}\tilde{R}%
_{kmq}^{\ell}$, we conclude
\begin{align*}
\tilde{\nabla}_{0}\tilde{R}_{ijk}^{\ell}  &  =\tilde{\Delta}\tilde{R}%
_{ijk}^{\ell}+2\tilde{B}_{ikj}^{\ell}-2\tilde{B}_{jki}^{\ell}\\
&  +\tilde{g}^{pq}\left(  \tilde{R}_{pji}^{m}-\tilde{R}_{pij}^{m}\right)
\left(  \tilde{R}_{kqm}^{\ell}-\tilde{R}_{kmq}^{\ell}\right)  +2\mu\tilde
{R}_{ijk}^{\ell}\\
&  =\tilde{\Delta}\tilde{R}_{ijk}^{\ell}+2\left(  \tilde{B}_{ikj}^{\ell
}-\tilde{B}_{jki}^{\ell}-\tilde{B}_{jik}^{\ell}+\tilde{B}_{ijk}^{\ell}\right)
+2\mu\tilde{R}_{ijk}^{\ell}.
\end{align*}

\end{proof}

\subsection[Space-time curvature]{Space-time curvature as a bilinear form}

We shall find it convenient to regard the curvature tensor as type $\left(
4,0\right)  .$ Since the space-time metric is degenerate, we lower indices as
follows:%
\begin{equation}
\tilde{R}_{ijk\ell}\doteqdot\left\{
\begin{tabular}
[c]{ll}%
$\bar{g}_{\ell p}\tilde{R}_{ijk}^{p}$ & if $\ell\geq1$\\
$-\bar{g}_{kp}\tilde{R}_{ij\ell}^{p}$ & if $\ell=0$ and $k\geq1$\\
$0$ & if $\ell=k=0.$%
\end{tabular}
\right.  . \label{Bilinear1}%
\end{equation}
We may now consider $\widetilde{\operatorname*{Rm}}$ to be a symmetric
quadratic form on $\Lambda^{2}T\widetilde{\mathcal{M}}$ by defining%
\begin{equation}
\widetilde{\operatorname*{Rm}}\left(  \tilde{S},\tilde{T}\right)
\doteqdot\sum_{i,j,k,\ell=0}^{n}\tilde{R}_{ijk\ell}\tilde{S}^{ij}\tilde
{T}^{\ell k}. \label{Bilinear2}%
\end{equation}
Note that this differs slightly from \cite{ChowChu1}, where $\widetilde
{\operatorname*{Rm}}$ was regarded as a tensor of type $\left(  2,2\right)  $.

Setting $\tilde{B}_{ijk\ell}\doteqdot-\tilde{g}^{pq}\tilde{R}_{pij}^{m}%
\tilde{R}_{kqm\ell}$, we restate Lemma \ref{st-riem-31-evol} in the form:

\begin{corollary}
\label{st-riem-evol-eqn-4-0}If $\mu$ is constant, then%
\begin{equation}
\tilde{\nabla}_{0}\tilde{R}_{ijk\ell}=\tilde{\Delta}\tilde{R}_{ijk\ell
}+2\left(  \tilde{B}_{ijk\ell}-\tilde{B}_{jik\ell}-\tilde{B}_{jki\ell}%
+\tilde{B}_{ikj\ell}\right)  +2\mu\tilde{R}_{ijk\ell}.
\label{st-riem-40-b-evol}%
\end{equation}

\end{corollary}

Recall that the degenerate metric $\tilde{g}$ induces an inner product and Lie
bracket on $\Lambda^{2}T^{\ast}\widetilde{\mathcal{M}}$ as follows:%
\begin{align}
\left\langle \tilde{S},\tilde{T}\right\rangle  &  =\tilde{g}^{ik}\tilde
{g}^{j\ell}\tilde{S}_{ij}\tilde{T}_{k\ell}\label{Winp}\\
\left[  \tilde{S},\tilde{T}\right]  _{ij}  &  =\tilde{g}^{k\ell}\left(
\tilde{S}_{ik}\tilde{T}_{\ell j}-\tilde{T}_{ik}\tilde{S}_{\ell j}\right)  .
\label{Wbrk}%
\end{align}
(Compare formulas (10) and (11) in \cite{ChowChu1}.) The structure constants
$C_{ij}^{ab,cd}$ defined by%
\[
C_{ij}^{ab,cd}dx^{i}\wedge dx^{j}\doteqdot\left[  dx^{a}\wedge dx^{b}%
,dx^{c}\wedge dx^{d}\right]
\]
for $0\leq a<b\leq n,$ $0\leq c<d\leq n$ and $0\leq i<j\leq n$ are then given
by%
\[
C_{ij}^{ab,cd}=\delta_{i}^{a}\delta_{j}^{d}\tilde{g}^{bc}-\delta_{i}^{c}%
\delta_{j}^{b}\tilde{g}^{ad}.
\]
In terms of the natural isomorphism between $\Lambda^{2}T^{\ast}%
\widetilde{\mathcal{M}}$ and $\Lambda^{2}T^{\ast}\mathcal{M}\oplus\Lambda
^{1}T^{\ast}\mathcal{M},$ formula (\ref{Wbrk}) corresponds to%
\[
\left[  X\oplus V,Y\oplus W\right]  =\left[  X,Y\right]  \oplus\left(
V\lrcorner Y-W\lrcorner X\right)  ,
\]
where $\lrcorner$ denotes the interior product. Analogous to \cite{HamPCO} and
the extension in \cite{ChowChu1}, we define a symmetric bilinear operator $\#$
on $\Lambda^{2}T^{\ast}\widetilde{\mathcal{M}}\otimes\Lambda^{2}T^{\ast
}\widetilde{\mathcal{M}}$ by%
\[
\left(  F\#G\right)  _{ijk\ell}\doteqdot F_{abcd}G_{pqrs}C_{ij}^{ab,pq}C_{\ell
k}^{cd,rs},
\]
and adopt the notational convenience $F^{\#}\doteqdot F\#F$. We also define
the square of an element in $\Lambda^{2}T^{\ast}\widetilde{\mathcal{M}}%
\otimes\Lambda^{2}T^{\ast}\widetilde{\mathcal{M}}$ by%
\[
F_{ijk\ell}^{2}\doteqdot\tilde{g}^{ad}\tilde{g}^{bc}F_{ijab}F_{cdk\ell}.
\]
With these definitions, following \cite{HamPCO} and \cite{ChowChu1}, we find
that (\ref{st-riem-40-b-evol}) takes the form:

\begin{lemma}
If $\mu$ is constant, then%
\begin{equation}
\tilde{\nabla}_{0}\tilde{R}_{ijk\ell}=\tilde{\Delta}\tilde{R}_{ijk\ell}%
+\tilde{R}_{ijk\ell}^{2}+\tilde{R}_{ijk\ell}^{\#}+2\mu\tilde{R}_{ijk\ell}.
\label{st-riem-sharo-evol}%
\end{equation}

\end{lemma}

We omit the long but straightforward computations.

\section[Hamilton's quadratic]{Hamilton's quadratic for the Ricci flow}

As was remarked in \cite{ChowChu1}, the results above give an explanation for
the surprising identities observed by Hamilton in \S 14 of \cite{HamForm}.
Here, we shall exhibit a correspondence between the machinery Hamilton uses to
prove his tensor inequality and the geometric structure of space-time, in
order to show that his quadratic and the assumptions made in its derivation
arise very naturally in the space-time context. A similar but less precise
correspondence appeared earlier in \cite{ChowChu1}.

Recall that Hamilton proved that for any $2$-form $U$ and $1$-form $W$ on a
complete solution $\left(  \mathcal{M}^{n},g\left(  t\right)  \right)  $ of
the Ricci flow with non-negative curvature operator, the\emph{ }quadratic%
\begin{equation}
Z=Z\left(  U,W\right)  \doteqdot M_{ij}W^{i}W^{j}+2P_{ijk}U^{ij}%
W^{k}+R_{ijk\ell}U^{ij}U^{\ell k} \label{Z}%
\end{equation}
is non-negative at all positive times, where $R_{ijk\ell}=g_{\ell m}%
R_{ijk}^{m}$,%
\begin{equation}
M_{ij}\doteqdot\Delta R_{ij}-\frac{1}{2}\nabla_{i}\nabla_{j}R+2R_{ipqj}%
R^{pq}-R_{ip}R_{j}^{p}+\frac{1}{2t}R_{ij}, \label{M}%
\end{equation}
and%
\begin{equation}
P_{ijk}\doteqdot\nabla_{i}R_{jk}-\nabla_{j}R_{ik}. \label{P}%
\end{equation}

We shall now relate Hamilton's proof to our construction. In \S 2 of
\cite{HamHarnack}, tensors of type $\left(  r,s\right)  $ on a Riemannian
manifold $\left(  \mathcal{M},g\right)  $ are regarded as $\operatorname*{GL}%
\left(  n,\mathbb{R}\right)  $-invariant maps from the linear frame bundle
$\operatorname*{GL}\left(  \mathcal{M}\right)  $ to $\mathbb{R}^{n^{r+s}}$.
For instance, if $P\in\mathcal{M}$ and $Y=\left(  Y_{1}|_{P},\ldots,Y_{n}%
|_{P}\right)  \in\operatorname*{GL}\left(  \mathcal{M}\right)  $ is given by
$Y_{a}=y_{a}^{i}\,\partial/\partial x^{i}$ in a chart $\left\{  x^{i}\right\}
$ at $P$, a $1$-form $\theta$ may be identified with the system of component
functions $\theta_{a}=\theta\left(  Y_{a}\right)  $ it induces on
$\operatorname*{GL}\left(  \mathcal{M}\right)  $. Regarding the Levi-Civita
connection of $\left(  \mathcal{M},g\right)  $ as a $\operatorname*{GL}\left(
n,\mathbb{R}\right)  $-invariant choice of horizontal subspace for each
$T_{Y}\operatorname*{GL}\left(  \mathcal{M}\right)  $, Hamilton takes
space-like derivatives by means of the unique horizontal lift $D_{a}$ of
$Y_{a}$ at $Y\in\operatorname*{GL}\left(  \mathcal{M}\right)  $. Hamilton then
identifies the vertical vector field $\vee_{a}^{b}$ with the differential of
the map $Y_{a}\mapsto y_{b}^{i}\left(  y^{-1}\right)  _{i}^{a}Y_{a}$;
namely\footnote{Hamilton writes $\nabla_{a}^{b}$ for what we denote $\vee
_{a}^{b}$.}
\[
\vee_{b}^{a}=y_{b}^{i}\frac{\partial}{\partial y_{a}^{i}}.
\]
Note that $\vee_{b}^{a}$ acts on a covariant tensor by
\begin{equation}
\vee_{b}^{a}T_{cd\dots z}=\delta_{c}^{a}T_{bd\dots z}+\delta_{d}^{a}T_{cb\dots
z}+\cdots+\delta_{z}^{a}T_{cd\dots b}. \label{vertical-VF}%
\end{equation}
For a solution $\left(  \mathcal{M},g\left(  t\right)  \right)  $ to the Ricci
flow on an interval $\mathcal{I}$, one considers the bundle
$\operatorname*{GL}\left(  \mathcal{M}\right)  \times\mathcal{I}$
$\rightarrow\widetilde{\mathcal{M}}=\mathcal{M}\times\mathcal{I}$ and the
sub-bundle of orthonormal frames $O\left(  \widetilde{\mathcal{M}}\right)
\doteqdot\cup_{t\in\mathcal{I}}\left(  O\left(  \mathcal{M},g\left(  t\right)
\right)  ,t\right)  \rightarrow\widetilde{\mathcal{M}}$. Hamilton takes
time-like derivatives by means of a vector field $D_{t}$ on
$\operatorname*{GL}\left(  \mathcal{M}\right)  \times\mathcal{I}$ defined by
\begin{equation}
D_{t}\doteqdot\frac{\partial}{\partial t}+R_{ab}g^{bc}\vee_{c}^{a}.
\label{Dt-definition}%
\end{equation}
$D_{t}$ is tangent to $O\left(  \widetilde{\mathcal{M}}\right)  $, because
\begin{equation}
D_{t}g_{ab}\equiv0. \label{D_t g_ab = 0}%
\end{equation}
The geometric structure of space-time reveals why this construction is
natural. Indeed, definition (\ref{Dt-definition}) corresponds to (\ref{C2}) in
the definition of the space-time connection $\tilde{\nabla}$ for the Ricci
flow without rescaling $\left(  \mu=0\right)  $, because%
\[
\tilde{\nabla}_{0}=\frac{\partial}{\partial\bar{t}}-\tilde{\Gamma}_{0j}%
^{k}\vee_{k}^{j}%
\]
when acting on covariant tensors. Property (\ref{D_t g_ab = 0}) corresponds to
the compatibility of $\tilde{\nabla}$ with the space-time metric $\tilde{g}$
(Lemma \ref{st-metric-parallel}).

Now suppose $\left(  \mathcal{M},g\left(  t\right)  \right)  $ is a
homothetically expanding soliton flowing along a gradient vector field $V$.
(See Definition \ref{ExpRcSoliton}, and recall that $\left(  \mathcal{M}%
,\bar{g}\left(  \bar{t}\right)  \right)  $ is then a steady Ricci soliton
flowing along $\bar{V}=e^{\bar{t}}V$.) In \S 3 of \cite{HamHarnack}, such a
solution is described by the equation%
\begin{equation}
D_{a}V_{b}=D_{b}V_{a}\triangleq R_{ab}+\frac{1}{2t}g_{ab}.
\label{expand-soliton-eqn}%
\end{equation}
Here and in what follows, we use the symbol $\triangleq$ to denote an identity
that holds for an expanding gradient soliton. By applying formula (\ref{P1})
to $\tilde{V}$, we note that condition (\ref{expand-soliton-eqn}) holds if and
only if for all $i,j\geq1$, one has%
\[
\tilde{\nabla}_{i}\tilde{V}^{j}=0.
\]

Hamilton next defines the quadratic $Z$ in terms of the tensors $M$, $P$, and
$\operatorname*{Rm}$. (Recall (\ref{Z})--(\ref{P}) and note that our sign
convention for the Riemann curvature tensor is opposite Hamilton's.) In
analogy with Theorem 2.2 of \cite{ChowChu1}, we apply Proposition
\ref{st-riem-formulas} with $\mu\equiv1/2$ to observe that these also
correspond to natural space-time objects:

\begin{lemma}
\label{RPM}Let $\left(  \mathcal{M},g\left(  t\right)  \right)  $ be a
solution of the Ricci flow. Set $\bar{t}=\ln t$ and $\bar{g}\left(  \bar
{t}\right)  =\frac{1}{t}g\left(  t\right)  $. Then for $i,j,k,\ell\geq1$, one
has
\begin{align}
R_{ijk\ell}  &  =e^{\bar{t}}\tilde{R}_{ijk\ell}\\
P_{\ell kj}  &  =\tilde{R}_{0jk\ell}=\tilde{R}_{k\ell0j}\label{R0jkl=Plkj}\\
M_{i\ell}  &  =e^{-\bar{t}}\tilde{R}_{i00\ell}. \label{Ri00l=Mil}%
\end{align}

\end{lemma}

\noindent Thus we arrive at the key observation that \emph{the LYH quadratic
may be identified with the space-time curvature tensor:}%
\[
Z=e^{\bar{t}}\sum_{i,j,k,\ell=0}^{n}\tilde{R}_{ijk\ell}\tilde{T}^{ij}\tilde
{T}^{\ell k},
\]
where the space-time contravariant $2$-tensor $\tilde{T}$ is defined in terms
of the natural isomorphism $\Lambda^{2}T^{\ast}\widetilde{\mathcal{M}}%
\cong\Lambda^{2}T^{\ast}\mathcal{M}\oplus\Lambda^{1}T^{\ast}\mathcal{M}$ by
$\tilde{T}\doteqdot U\oplus\left(  e^{-\bar{t}}W/2\right)  $. In components
$i,j\geq1$,%
\begin{align}
\tilde{T}^{ij}  &  =U^{ij}\label{Tij-def}\\
\tilde{T}^{0j}  &  =-\tilde{T}^{j0}=\frac{1}{2}e^{-\bar{t}}W^{j}=\frac{1}%
{2t}W^{j}. \label{T0j-def}%
\end{align}
(See also Corollary 2.3 of \cite{ChowChu1}; a key difference from that paper
is that taking $\mu=\frac{1}{2}$ accounts for the term $\frac{1}{2t}R_{i\ell}$
in $M_{i\ell}$.)

Differentiating the expanding gradient soliton equation
(\ref{expand-soliton-eqn}), Hamilton obtains the following two relations:%
\begin{align}
P_{abc}+R_{abc}^{d}V_{d}  &  \triangleq0\label{Pabc-sol-eqn}\\
M_{ab}+P_{cab}V^{c}  &  \triangleq0. \label{Mab-sol-eqn}%
\end{align}
Together, these equations prove the Li--Yau--Hamilton inequality is sharp.
Indeed, if $W$ is arbitrary and one sets $U_{ab}=\frac{1}{2}\left(  V_{a}%
W_{b}-V_{b}W_{a}\right)  $, a straightforward computation gives $Z\left(
U,W\right)  \triangleq0$. This fact can be interpreted using the result of
Corollary \ref{Vparallel} that%
\begin{equation}
\tilde{R}_{ijk}^{\ell}\tilde{V}^{k}\triangleq0 \label{st-Riem-vanishes-V}%
\end{equation}
holds for all $i,j,\ell\geq0$:

\begin{lemma}
The identities (\ref{Pabc-sol-eqn}) and (\ref{Mab-sol-eqn}) are equivalent to
the fact that the space-time Riemannian curvature tensor $\widetilde
{\operatorname*{Rm}}$ vanishes in the direction of the parallel space-time
vector field $e^{\bar{t}/2}\tilde{V}$ when $\mu=1/2$.
\end{lemma}

\begin{proof}
If $i,j,\ell\geq1$, Lemma \ref{RPM} implies that
\[
\tilde{R}_{ijk\ell}\tilde{V}^{k}=-e^{\bar{t}/2}\left(  P_{ij\ell}+R_{ij\ell
k}V^{k}\right)
\]
and
\[
\tilde{R}_{0jk\ell}\tilde{V}^{k}=-e^{3\bar{t}/2}\left(  M_{j\ell}+P_{k\ell
j}V^{k}\right)  .
\]
Since $\tilde{R}_{ijk}^{0}=0$ for all $i,j,k$, it is clear that
(\ref{st-Riem-vanishes-V}) holds if and only if both (\ref{Pabc-sol-eqn}) and
(\ref{Mab-sol-eqn}) do.
\end{proof}

The evolution equations satisfied by the coefficients of Hamilton's quadratic
are derived in Lemmas 4.2, 4.3, and 4.4 of \cite{HamHarnack}. Written in
Hamilton's notation, they are%
\begin{equation}
(D_{t}-\triangle)R_{abcd}=2(B_{abcd}-B_{abdc}+B_{acbd}-B_{adbc}),
\label{R-evolve}%
\end{equation}%
\begin{align}
(D_{t}-\triangle)P_{abc}  &  =-2R_{de}D_{d}R_{abce}\label{P-evolve}\\
&  +2\left(  R_{adbe}P_{dec}+R_{adce}P_{dbe}+R_{bdce}P_{ade}\right)
,\nonumber
\end{align}
and%
\begin{align}
(D_{t}-\triangle)M_{ab}  &  =2R_{cd}\left(  D_{c}P_{dab}+D_{c}P_{dba}\right)
+2R_{acbd}M_{cd}\label{M-evolve}\\
&  +2P_{acd}P_{bcd}-4P_{acd}P_{bdc}+2R_{cd}R_{ce}R_{adbe}-{\frac{1}{2t^{2}}%
}R_{ab},\nonumber
\end{align}
where $B_{abcd}=R_{aebf}R_{cedf}$. In \S \ref{Evolution}, we prove the following:

\begin{proposition}
\label{equivalence}The evolution equations (\ref{R-evolve}), (\ref{P-evolve}),
and (\ref{M-evolve}) are equivalent to the evolution equation%
\begin{equation}
\tilde{\nabla}_{0}\tilde{R}_{ijk\ell}=\tilde{\Delta}\tilde{R}_{ijk\ell
}+2\left(  \tilde{B}_{ijk\ell}-\tilde{B}_{jik\ell}-\tilde{B}_{jki\ell}%
+\tilde{B}_{ikj\ell}\right)  +\tilde{R}_{ijk\ell}
\label{st-riem-evol-eqn-40-II}%
\end{equation}
satisfied by $\widetilde{\operatorname*{Rm}}$ when $\mu=1/2$.
\end{proposition}

In computing the evolution of the quadratic $Z,$ Hamilton makes the following
assumptions on the $2$-form $U$ and the $1$-form $W$ at a given point:%
\begin{align}
\left(  D_{t}-\Delta\right)  W_{a}  &  =\frac{1}{t}W_{a}\tag{A1}\label{A1}\\
\left(  D_{t}-\Delta\right)  U_{ab}  &  =0\tag{A2}\label{A2}\\
D_{a}W_{b}  &  =0\tag{A3}\label{A3}\\
D_{a}U_{bc}  &  =\frac{1}{2}\left(  R_{ab}W_{c}-R_{ac}W_{b}\right)  +\frac
{1}{4t}\left(  g_{ab}W_{c}-g_{ac}W_{b}\right)  . \tag{A4}\label{A4}%
\end{align}
(See the hypotheses in Theorem 4.1 of \cite{HamHarnack}, and note that
equation (\ref{A4}) is motivated by the fact that it holds on a soliton if
(\ref{A3}) holds and $U=V\wedge W$.) We shall now demonstrate that \emph{the
four assumptions above are also very natural from the space-time perspective.
}Indeed, equations (\ref{A1})--(\ref{A4}) hold at a point in space-time if and
only if $e^{\bar{t}}\tilde{T}^{ij}$ satisfies the heat equation and is
parallel in space-like directions at that point:

\begin{lemma}
If $\mu=1/2$, assumptions (\ref{A1})-(\ref{A4}) are equivalent to
\begin{align}
\left(  \tilde{\nabla}_{0}-\tilde{\Delta}+1\right)  \tilde{T}^{ij}  &
=0\label{heat-Tij}\\
\tilde{\nabla}_{k}\tilde{T}^{ij}  &  =0, \label{covar-Tij}%
\end{align}
for all $i,j\geq0$ and $k\geq1$.
\end{lemma}

\begin{proof}
For $i,j,k\geq1$, we use (\ref{C1}), (\ref{C2}), (\ref{Tij-def}),
(\ref{T0j-def}), and the fact that $\bar{R}_{k}^{i}=tR_{k}^{i}$ to compute%
\begin{align*}
\tilde{\nabla}_{k}\tilde{T}^{ij}  &  =\nabla_{k}U^{ij}+\tilde{\Gamma}_{k0}%
^{i}\tilde{T}^{0j}+\tilde{\Gamma}_{k0}^{j}\tilde{T}^{i0}\\
&  =\nabla_{k}U^{ij}-\left(  tR_{k}^{i}+\frac{1}{2}\delta_{k}^{i}\right)
\frac{1}{2t}W^{j}+\left(  tR_{k}^{j}+\frac{1}{2}\delta_{k}^{j}\right)
\frac{1}{2t}W^{i}.
\end{align*}
Hence (\ref{covar-Tij}) is valid for all $i,j,k\geq1$ if and only if
(\ref{A4}) holds. For $i=0$ but $j,k\geq1$, we have%
\[
\tilde{\nabla}_{k}\tilde{T}^{0j}=\frac{1}{2t}\nabla_{k}W^{j}.
\]
So (\ref{covar-Tij}) is valid for $ij=0$ and all $k\geq1$ if and only if
(\ref{A3}) holds. Similarly, since $\tilde{\nabla}_{q}\tilde{T}^{0j}=\frac
{1}{2t}\nabla_{q}W^{j}$ and $\tilde{\Delta}\tilde{T}^{0j}=\tilde{g}^{pq}%
\tilde{\nabla}_{p}\tilde{\nabla}_{q}\tilde{T}^{0j}=\frac{1}{2}\Delta W^{j}$,
we compute that%
\begin{align*}
\left(  \tilde{\nabla}_{0}-\tilde{\Delta}+1\right)  \tilde{T}^{0j}  &
=\frac{\partial}{\partial\bar{t}}\left(  \frac{1}{2t}W^{j}\right)
+\tilde{\Gamma}_{00}^{0}\tilde{T}^{0j}+\tilde{\Gamma}_{0p}^{j}\tilde{T}^{0p}\\
&  -\frac{1}{2}\Delta W^{j}+\frac{1}{2t}W^{j}\\
&  =\frac{1}{2}\left(  \frac{\partial}{\partial t}W^{j}-R_{p}^{j}W^{p}-\Delta
W^{j}-\frac{1}{t}W^{j}\right)  .
\end{align*}
It follows easily that (\ref{heat-Tij}) is valid for $ij=0$ if and only if
(\ref{A1}) holds. Finally, we use (\ref{C3}) to calculate%
\begin{align*}
\left(  \tilde{\nabla}_{0}+1\right)  \tilde{T}^{ij}  &  =\frac{\partial
}{\partial\bar{t}}U^{ij}+\tilde{\Gamma}_{0p}^{i}U^{pj}+\tilde{\Gamma}_{00}%
^{i}\tilde{T}^{0j}+\tilde{\Gamma}_{0p}^{j}U^{ip}+\tilde{\Gamma}_{00}^{j}%
\tilde{T}^{i0}+U^{ij}\\
&  =t\left(  \frac{\partial}{\partial t}U^{ij}-R_{p}^{i}U^{pj}-R_{p}^{j}%
U^{ip}\right)  -\frac{1}{4}W^{j}\nabla^{i}R+\frac{1}{4}W^{i}\nabla^{j}R.
\end{align*}
Then noting that for $i,j\geq1$,%
\[
\tilde{\nabla}_{q}\tilde{T}^{ij}=\nabla_{q}U^{ij}+\frac{1}{2}\left(  R_{q}%
^{j}W^{i}-R_{q}^{i}W^{j}\right)  +\frac{1}{4t}\left(  \delta_{q}^{j}%
W^{i}-\delta_{q}^{i}W^{j}\right)  ,
\]
we compute%
\begin{align*}
\tilde{\Delta}\tilde{T}^{ij}  &  =\tilde{g}^{pq}\tilde{\nabla}_{p}%
\tilde{\nabla}_{q}\tilde{T}^{ij}\\
&  =tg^{pq}\left[
\begin{array}
[c]{c}%
\nabla_{p}\nabla_{q}U^{ij}+\frac{1}{2}\nabla_{p}\left(  R_{q}^{j}W^{i}%
-R_{q}^{i}W^{j}\right) \\
+\frac{1}{4t}\nabla_{p}\left(  \delta_{q}^{j}W^{i}-\delta_{q}^{i}W^{j}\right)
+\tilde{\Gamma}_{p0}^{i}\tilde{\nabla}_{q}\tilde{T}^{0j}+\tilde{\Gamma}%
_{p0}^{j}\tilde{\nabla}_{q}\tilde{T}^{i0}%
\end{array}
\right] \\
&  =t\left[  \Delta U^{ij}+R_{p}^{j}\nabla^{p}W^{i}-R_{p}^{i}\nabla^{p}%
W^{j}+\frac{1}{4}\left(  W^{i}\nabla^{j}R-W^{j}\nabla^{i}R\right)  \right] \\
&  +\frac{1}{2}\left(  \nabla^{j}W^{i}-\nabla^{i}W^{j}\right)
\end{align*}
and collect terms to obtain%
\begin{align*}
\left(  \tilde{\nabla}_{0}-\tilde{\Delta}+1\right)  \tilde{T}^{ij}  &
=t\left[  \frac{\partial}{\partial t}U^{ij}-\Delta U^{ij}-R_{p}^{i}%
U^{pj}-R_{p}^{j}U^{ip}\right] \\
&  +\left(  tR^{ip}+\frac{1}{2}g^{ip}\right)  \nabla_{p}W^{j}-\left(
tR^{jp}+\frac{1}{2}g^{jp}\right)  \nabla_{p}W^{i}.
\end{align*}
So if (\ref{A3}) holds, then (\ref{heat-Tij}) is valid for $i,j\geq1$ if and
only if (\ref{A2}) holds.
\end{proof}

\section{A generalized tensor maximum principle}

In order to utilize space-time methods fully in investigating potential
Li--Yau--Hamilton quadratics for the Ricci flow, one needs a version of the
parabolic maximum principle for equations such as (\ref{st-riem-sharo-evol})
and (\ref{gen-Rm-evolve-eq1}) Accordingly, we now derive a generalization of
the tensor maximum principle originally proved in \cite{Ham3man}. We begin
with the observation that any smooth family $\left\{  g\left(  t\right)
:0\leq t<\Omega\right\}  $ of Riemannian metrics on $\mathcal{M}^{n}$ induces
a nondegenerate metric $\hat{g}$ on $\mathcal{M}\times\lbrack0,\Omega)$ given
in coordinates $\left(  \partial/\partial t=x^{0},\,x^{1},\dots,x^{n}\right)
$ by%
\[
\hat{g}_{ij}=\left\{
\begin{array}
[c]{cl}%
g_{ij} & \text{if }1\leq i,j\leq n\\
\delta_{ij} & \text{if }i=0\text{ or }j=0.
\end{array}
\right.
\]
We denote the Levi-Civita connection of $\hat{g}$ by $\hat{\nabla}$.

\begin{proposition}
\label{STMP}Let $g\left(  t\right)  $ be a smooth $1$-parameter family of
complete metrics on $\mathcal{M}^{n}$, indexed by $t\in\lbrack0,\Omega)$. Let
$\widetilde{\mathcal{M}}\doteqdot\mathcal{M}\times\lbrack0,\Omega)$ and let
$\tilde{g}$ be the degenerate metric defined on $T^{\ast}\widetilde
{\mathcal{M}}$ by%
\[
\tilde{g}^{ij}\doteqdot\left\{
\begin{tabular}
[c]{ll}%
$g^{ij}$ & if $\text{if }1\leq i,j\leq n$\\
$0$ & if $i=0$ or $j=0.$%
\end{tabular}
\right.
\]
Let $\tilde{\nabla}$ be a compatible connection $(\tilde{\nabla}_{i}\tilde
{g}^{jk}\equiv0)$, and let $\mathcal{Q}$ denote the space of symmetric
bilinear forms on a tensor bundle $\mathcal{X}$ over $\widetilde{\mathcal{M}}%
$. Suppose $Q\in\mathcal{Q}$\ is a solution of the reaction-diffusion
equation
\begin{equation}
\tilde{\nabla}_{0}Q=\tilde{\Delta}Q+\Phi\left(  Q\right)  , \label{RDE}%
\end{equation}
where $\Phi:\mathcal{Q}\rightarrow\mathcal{Q}$\ is a (possibly nonlinear)
locally Lipschitz map which satisfies the \emph{null eigenvector condition
}that $\Phi\left(  P\right)  \left(  X,X\right)  \geq0$ at any point where
$P\left(  X,\cdot\right)  $ vanishes for $P\in\mathcal{Q}$ and $X\in
\mathcal{X}$. Assume $\left|  \tilde{\nabla}-\hat{\nabla}\right|  _{\hat{g}}$,
$\left|  \hat{\nabla}\left(  \tilde{\nabla}-\hat{\nabla}\right)  \right|
_{\hat{g}}$, and the Lipschitz constant for $\Phi$ are bounded on any subset
$\mathcal{M}\times\left[  0,\eta\right]  \subset\widetilde{\mathcal{M}}$. If
$\mathcal{M}$ is not compact, assume also that there exists $\rho
:\mathcal{M}\rightarrow\lbrack1,\infty)$ with $\rho^{-1}\left(  \left[
1,s\right]  \right)  $ compact for every $s\in\lbrack1,\infty)$ and such that
$\left|  \nabla\rho\right|  _{g}$ and $\left|  \Delta\rho\right|  $ are
bounded on $\mathcal{M}\times\left[  0,\eta\right]  $. If $Q\geq0$ on
$\mathcal{M}\times\left\{  0\right\}  $. then $Q\geq0$ on $\widetilde
{\mathcal{M}}$.
\end{proposition}

\begin{proof}
The metric $\hat{g}$ induces an inner product on $\mathcal{X}$ in the usual
way; we shall abuse notation by writing $\hat{g}\left(  X,Y\right)  $ for
$X,Y\in\mathcal{X}$. If $\mathcal{M}$ is compact, take $\rho\equiv1$, and
otherwise let $\rho:\widetilde{\mathcal{M}}\rightarrow\lbrack1,\infty)$ be the
function in the statement of the theorem. (By \cite{GreeneWu} and Lemma 5.1 of
\cite{HamHarnack}, such a function always exists if the time derivatives
$\partial g/\partial t$ of $g$ and $\partial\Gamma/\partial t$ of the
Levi-Civita connection of $g$ are bounded, and if $\mathcal{M}$ has positive
sectional curvature.)

By considering translates in time, it will suffice to prove there is $\eta>0$
such that for every $\varepsilon>0$, the quadratic form $\hat{Q}$ is strictly
positive on $\mathcal{M}\times\left[  0,\eta\right]  $, where
\[
\hat{Q}\left(  x,t\right)  \doteqdot Q\left(  x,t\right)  +\varepsilon\left(
\eta+t\right)  \rho\left(  x\right)  \hat{g}\left(  x,t\right)  .
\]
Suppose $\hat{Q}$ does not remain strictly positive, and let $t_{0}\in\left[
0,\eta\right]  $ denote the infimum of all $t$ such that $\left.  \hat
{Q}\left(  Y,Y\right)  \right|  _{\left(  x,t\right)  }=0$ for some
$Y\in\mathcal{X}$ and $x\in\mathcal{M}$ with $\left|  Y\right|  _{\hat{g}}=1$
at $\left(  x,t\right)  $. We claim $t_{0}>0$. If not, there will be a
sequence of compact sets $\mathcal{K}_{j}$ exhausting $\mathcal{M}$, points
$x_{j}\in\mathcal{K}_{j}\backslash\mathcal{K}_{j-1}$, and times $t_{j}\searrow
t_{0}=0$ such that the first zero of $\hat{Q}$ on $\mathcal{K}_{j}%
\times\left[  0,\eta\right]  $ occurs at $\left(  x_{j},t_{j}\right)  $. Since
$Q\geq0$ on $\mathcal{M}\times\left\{  0\right\}  $ and $\rho\left(
x_{j}\right)  \rightarrow\infty$ if $\mathcal{M}$ is not compact, this is impossible.

By the null eigenvector assumption,%
\[
\left.  \Phi(\hat{Q})\left(  Y,Y\right)  \right|  _{\left(  x,t_{0}\right)
}\geq0.
\]
Define a tensor field $X$ in a space-time neighborhood $\mathcal{O}$ of
$\left(  x,t_{0}\right)  $ by taking $X=Y$ at $\left(  x,t_{0}\right)  $ and
extending $X$ by parallel transport along radial geodesics with respect to the
connection $\tilde{\nabla}$. (It suffices to extend $X$ first radially along
all $\tilde{\nabla}$-geodesics which start tangent to the hypersurface
$\mathcal{M}\times\left\{  t_{0}\right\}  $, and then along any curve with
tangent $\partial/\partial t$ at $\left(  x,t_{0}\right)  $.) Notice that all
symmetric space-like second covariant derivatives of $X$ vanish at $\left(
x,t_{0}\right)  $. (Compare \S 4 of \cite{HamPCO}.) Indeed, with respect to a
$\hat{g}$-orthonormal frame $\left\{  e_{0}=\partial/\partial t,\,e_{1}%
,\dots,e_{n}\right\}  $, one observes that for $i=1,\dots,n$,%
\[
\left.  \tilde{\nabla}_{e_{i}}\tilde{\nabla}_{e_{i}}X\right|  _{\left(
x,t_{0}\right)  }=\tilde{\nabla}_{e_{i}}\left(  \tilde{\nabla}_{e_{i}%
}X\right)  -\tilde{\nabla}_{\tilde{\nabla}_{e_{i}}e_{i}}X=0-0.
\]
Hence for any $P\in\mathcal{Q}$, we compute at $\left(  x,t_{0}\right)  $ that%
\begin{align*}
\tilde{\Delta}\left(  P\left(  X,X\right)  \right)   &  =\tilde{g}^{ij}\left[
\begin{array}
[c]{c}%
(\tilde{\nabla}_{i}\tilde{\nabla}_{j}P)(X,X)+4(\tilde{\nabla}_{i}%
P)(X,\tilde{\nabla}_{j}X)\\
+2P(\tilde{\nabla}_{i}X,\tilde{\nabla}_{j}X)+2P(X,\tilde{\nabla}_{i}%
\tilde{\nabla}_{j}X)
\end{array}
\right] \\
&  =(\tilde{\Delta}P)\left(  X,X\right)  .
\end{align*}

Now consider the function $F$ defined in $\mathcal{O}$ by%
\[
F\left(  y,t\right)  =\left.  \hat{Q}\left(  X,X\right)  \right|  _{\left(
y,t\right)  }.
\]
Even though $\hat{g}$ may not be compatible with the connection $\tilde
{\nabla}$, we still have $\left|  X\right|  _{\hat{g}}\geq1/2$ in a possibly
smaller neighborhood $\mathcal{O}^{\prime}\subseteq\mathcal{O}$. Hence $F$
attains its minimum in $\mathcal{O}^{\prime}\cap\mathcal{M}\times\left[
0,t_{0}\right]  $ at $\left(  x,t_{0}\right)  $, where we therefore have%
\[
0\geq\frac{\partial}{\partial t}F=\left(  \tilde{\nabla}_{0}\hat{Q}\right)
\left(  X,X\right)  ,
\]
and
\[
0=\frac{\partial}{\partial x^{i}}F=\left(  \tilde{\nabla}_{i}\hat{Q}\right)
\left(  X,X\right)
\]
for $i=1,\dots,n$, and%
\[
0\leq\tilde{g}^{ij}\left(  \frac{\partial^{2}F}{\partial x^{i}\partial x^{j}%
}-\tilde{\Gamma}_{ij}^{k}\frac{\partial F}{\partial x^{k}}\right)
=\tilde{\Delta}F.
\]

To finish the proof, observe that there are constants $C_{1}$ and $C_{2}$
depending only on the bounds for $\left|  \tilde{\nabla}-\hat{\nabla}\right|
_{\hat{g}}$ and $\left|  \hat{\nabla}\left(  \tilde{\nabla}-\hat{\nabla
}\right)  \right|  _{\hat{g}}$on $\mathcal{M}\times\left[  0,\eta\right]  $
such that%
\[
2\left.  \left(  \tilde{\nabla}_{0}\hat{g}\right)  \left(  X,X\right)
\right|  _{\left(  x,t_{0}\right)  }\geq-C_{1}\left|  X\right|  _{\hat{g}}%
^{2}=-C_{1}%
\]
and%
\[
2\left.  \left(  \tilde{\Delta}\hat{g}\right)  \left(  X,X\right)  \right|
_{\left(  x,t_{0}\right)  }\leq C_{2}\left|  X\right|  _{\hat{g}}^{2}=C_{2}.
\]
There is $C_{3}$ depending only on the Lipschitz constant of $\Phi$ on
$\mathcal{M}\times\left[  0,\eta\right]  $ such that%
\[
-\left.  \Phi\left(  Q\right)  \left(  X,X\right)  \right|  _{\left(
x,t_{0}\right)  }\leq\Phi(\hat{Q})\left(  X,X\right)  -\Phi\left(  Q\right)
\left(  X,X\right)  \leq C_{3}\varepsilon\eta\left|  X\right|  _{\hat{g}}%
^{4}=\varepsilon\eta C_{3},
\]
and there is $C_{4}$ depending only on the bounds for $\left|  \Delta
\rho\right|  $, $\left|  \nabla\rho\right|  _{g}$, and $\left|  \tilde{\nabla
}-\hat{\nabla}\right|  _{\hat{g}}$ such that%
\begin{align*}
\left(  \tilde{\Delta}F\right)  \left(  x,t_{0}\right)   &  =\left(
\tilde{\Delta}Q\right)  \left(  X,X\right)  +\varepsilon\left(  \eta
+t_{0}\right)  \left[
\begin{array}
[c]{c}%
\rho\left(  \tilde{\Delta}\hat{g}\right)  \left(  X,X\right) \\
+2\left(  \tilde{\nabla}_{\tilde{\nabla}\rho}\hat{g}\right)  \left(
X,X\right) \\
+\left(  \Delta\rho\right)  \hat{g}\left(  X,X\right)
\end{array}
\right] \\
&  \leq\left(  \tilde{\Delta}Q\right)  \left(  X,X\right)  +\varepsilon
\eta\left(  \rho C_{2}+C_{4}\right)  .
\end{align*}
Combining these estimates with equation (\ref{RDE}), we conclude that at
$\left(  x,t_{0}\right)  $,%
\begin{align*}
0  &  \geq\frac{\partial}{\partial t}F\\
&  =\left(  \tilde{\nabla}_{0}Q\right)  \left(  X,X\right)  +\varepsilon
\left(  \eta+t_{0}\right)  \rho\left(  \tilde{\nabla}_{0}\hat{g}\right)
\left(  X,X\right)  +\varepsilon\rho\hat{g}\left(  X,X\right) \\
&  =\left(  \tilde{\Delta}Q\right)  \left(  X,X\right)  +\Phi\left(  Q\right)
\left(  X,X\right)  +\varepsilon\left(  \eta+t_{0}\right)  \rho\left(
\tilde{\nabla}_{0}\hat{g}\right)  \left(  X,X\right)  +\varepsilon\rho\\
&  \geq\tilde{\Delta}F-\varepsilon\eta\left(  \rho C_{2}+C_{4}\right)
-\varepsilon\eta C_{3}-\varepsilon\eta\rho C_{1}+\varepsilon\rho\\
&  \geq\varepsilon\left[  \rho\left(  1-\eta\left(  C_{1}+C_{2}\right)
\right)  -\eta\left(  C_{3}+C_{4}\right)  \right]  .
\end{align*}
Because $\rho\geq1$ and the constants $C_{i}$ cannot increase if $\eta$
decreases, choosing $\eta>0$ sufficiently small gives a contradiction. So
$\hat{Q}$ remains strictly positive on $\mathcal{M}\times\left[
0,\eta\right]  $, and the result follows by letting $\varepsilon\searrow0$.
\end{proof}

\section[Evolution equations]{Evolution equations relating to Hamilton's
quadratic\label{Evolution}}

This section is devoted to the following:

\begin{proof}
[Proof of Proposition \ref{equivalence}]Assume $\mu\equiv1/2$, and denote the
RHS of (\ref{st-riem-40-b-evol}) by $\digamma_{ijk\ell}$. Lemma \ref{RPM}
implies that the following identities are valid for all $i,j,k,\ell\geq1$:%
\begin{align}
\tilde{B}_{ijk0}  &  \doteqdot-\tilde{g}^{pq}\tilde{R}_{pij}^{r}\tilde
{R}_{kqr0}=tg^{pq}R_{pij}^{r}P_{qkr}\label{B-ident1}\\
\tilde{B}_{i00\ell}  &  \doteqdot-\tilde{g}^{pq}\tilde{R}_{pi0}^{r}\tilde
{R}_{0qr\ell}=-t^{2}g^{pq}g^{rs}P_{ipr}P_{\ell sq}\label{B-ident2}\\
\tilde{B}_{i0\ell0}  &  \doteqdot-\tilde{g}^{pq}\tilde{R}_{pi0}^{r}\tilde
{R}_{\ell qr0}=t^{2}g^{pq}g^{rs}P_{ipr}P_{q\ell s}\label{B-ident3}\\
\tilde{B}_{i\ell00}  &  \doteqdot-\tilde{g}^{pq}\tilde{R}_{pi\ell}^{r}%
\tilde{R}_{0qr0}=-t^{2}g^{pq}R_{pi\ell}^{r}M_{qr}. \label{B-ident4}%
\end{align}
Thus for $i,\ell\geq1$, we have%
\begin{align*}
\digamma_{i00\ell}  &  \doteqdot\tilde{\Delta}\tilde{R}_{i00\ell}+2\left(
2\tilde{B}_{i00\ell}-\tilde{B}_{0i0\ell}-\tilde{B}_{00i\ell}\right)
+\tilde{R}_{i00\ell}\\
&  =\tilde{\Delta}\tilde{R}_{i00\ell}+2t^{2}g^{pq}g^{rs}\left[  R_{ipr\ell
}M_{qs}-P_{ipr}\left(  2P_{\ell sq}+P_{q\ell s}\right)  \right]  +tM_{i\ell}.
\end{align*}
On the other hand, we use (\ref{C2})--(\ref{C4}) with Lemma \ref{RPM} to
compute directly that%
\begin{align*}
\tilde{\nabla}_{0}\tilde{R}_{i00\ell}  &  =\frac{\partial}{\partial\bar{t}%
}\tilde{R}_{i00\ell}-\tilde{\Gamma}_{0i}^{p}\tilde{R}_{p00\ell}-\tilde{\Gamma
}_{0\ell}^{p}\tilde{R}_{i00p}-\tilde{\Gamma}_{00}^{p}\left(  \tilde
{R}_{ip0\ell}+\tilde{R}_{i0p\ell}\right) \\
&  =t\frac{\partial}{\partial t}\left(  tM_{i\ell}\right)  +\left(  tR_{i}%
^{p}+\mu\delta_{i}^{p}\right)  \left(  tM_{p\ell}\right)  +\left(  tR_{\ell
}^{p}+\mu\delta_{\ell}^{p}\right)  \left(  tM_{ip}\right) \\
&  +\left(  \frac{1}{2}t^{2}\nabla^{p}R\right)  \left(  P_{pi\ell}-P_{\ell
pi}\right)  +2\mu\left(  tM_{i\ell}\right) \\
&  =t^{2}\left[  \frac{\partial}{\partial t}M_{i\ell}+R_{i}^{p}M_{p\ell
}+R_{\ell}^{p}M_{ip}+\frac{1}{2}\left(  \nabla^{p}R\right)  \left(  P_{pi\ell
}+P_{p\ell i}\right)  \right] \\
&  +3tM_{i\ell}.
\end{align*}
In the same way, we compute for $q,m\leq1$ that%
\[
\tilde{\nabla}_{q}\tilde{R}_{i00\ell}=t\left(  \nabla_{q}M_{i\ell}+R_{q}%
^{m}P_{mi\ell}+R_{q}^{m}P_{m\ell i}\right)  +\mu\left(  P_{qi\ell}+P_{q\ell
i}\right)
\]
and%
\[
\tilde{\nabla}_{q}\tilde{R}_{im0\ell}=\nabla_{q}P_{mi\ell}+R_{q}^{r}%
R_{imr\ell}+\frac{\mu}{t}R_{imq\ell}.
\]
Then by using the divergence identity%
\[
\nabla^{q}P_{qi\ell}=M_{i\ell}-R^{pq}R_{ipq\ell}-\frac{1}{2t}R_{i\ell},
\]
we can write%
\begin{align*}
\tilde{\Delta}\tilde{R}_{i00\ell}  &  =\bar{g}^{pq}\tilde{\nabla}_{p}%
\tilde{\nabla}_{q}\tilde{R}_{i00\ell}\\
&  =t^{2}\nabla^{q}\left(  \nabla_{q}M_{i\ell}+R_{q}^{m}P_{mi\ell}+R_{q}%
^{m}P_{m\ell i}\right)  +t\mu\nabla^{q}\left(  P_{qi\ell}+P_{q\ell i}\right)
\\
&  +tg^{pq}\left(  tR_{p}^{m}+\mu\delta_{p}^{m}\right)  \left(  \nabla
_{q}P_{mi\ell}+R_{q}^{r}R_{imr\ell}+\frac{\mu}{t}R_{imq\ell}\right) \\
&  +tg^{pq}\left(  tR_{p}^{m}+\mu\delta_{p}^{m}\right)  \left(  \nabla
_{q}P_{m\ell i}+R_{q}^{r}R_{irm\ell}+\frac{\mu}{t}R_{iqm\ell}\right) \\
&  =t^{2}\Delta M_{i\ell}+\frac{t^{2}}{2}\left(  \nabla^{q}R\right)  \left(
P_{qi\ell}+P_{q\ell i}\right)  +2t^{2}R^{pq}\nabla_{p}\left(  P_{qi\ell
}+P_{q\ell i}\right) \\
&  +2t^{2}R_{m}^{p}R^{mq}R_{ipq\ell}+2tM_{i\ell}-\frac{1}{2}R_{i\ell}.
\end{align*}
Cancelling terms yields%
\begin{align*}
\tilde{\nabla}_{0}\tilde{R}_{i00\ell}  &  =\tilde{\Delta}\tilde{R}_{i00\ell
}+t^{2}\left(  \frac{\partial}{\partial t}-\Delta\right)  M_{i\ell}%
-2t^{2}R^{pq}\nabla_{p}\left(  P_{qi\ell}+P_{q\ell i}\right) \\
&  +t^{2}\left(  R_{i}^{p}M_{p\ell}+R_{\ell}^{p}M_{ip}-2R_{m}^{p}%
R^{mq}R_{ipq\ell}\right)  +tM_{i\ell}+\frac{1}{2}R_{i\ell}.
\end{align*}
Recalling (\ref{vertical-VF}) and (\ref{Dt-definition}), we conclude that the
special case $\tilde{\nabla}_{0}\tilde{R}_{i00\ell}=\digamma_{i00\ell}$ of
equation (\ref{st-riem-40-b-evol}) holds if and only if%
\begin{align*}
D_{t}M_{i\ell}  &  \doteqdot\frac{\partial}{\partial t}M_{i\ell}+R_{p}^{q}%
\vee_{q}^{p}M_{i\ell}=\frac{\partial}{\partial t}M_{i\ell}+R_{i}^{p}M_{p\ell
}+R_{l}^{p}M_{ip}\\
&  =\Delta M_{i\ell}+2R^{pq}\nabla_{p}\left(  P_{qi\ell}+P_{q\ell i}\right)
+2R_{m}^{p}R^{mq}R_{ipq\ell}\\
&  +2g^{pq}g^{rs}\left[  R_{ipr\ell}M_{qs}-P_{ipr}\left(  2P_{\ell
sq}+P_{q\ell s}\right)  \right]  -\frac{1}{2t^{2}}R_{i\ell},
\end{align*}
hence if and only if equation (\ref{M-evolve}) holds.

Now if $i,j,k\geq1$, identities (\ref{B-ident1})--(\ref{B-ident4}) let us
write%
\begin{align*}
\digamma_{ijk0}  &  =\tilde{\Delta}\tilde{R}_{ijk0}+2tg^{pq}\left(
R_{pij}^{r}P_{qkr}-R_{pji}^{r}P_{qkr}-R_{pjk}^{r}P_{qir}+R_{pik}^{r}%
P_{qjr}\right) \\
&  +P_{ijk}\\
&  =\tilde{\Delta}\tilde{R}_{ijk0}+2tg^{pq}\left(  R_{pij}^{r}P_{qrk}%
-R_{pjk}^{r}P_{qir}+R_{pik}^{r}P_{qjr}\right)  +P_{ijk}.
\end{align*}
On the other hand, (\ref{C2})--(\ref{C4}) and Lemma \ref{RPM} imply that%
\begin{align*}
\tilde{\nabla}_{0}\tilde{R}_{ijk0}  &  =\frac{\partial}{\partial\bar{t}}%
\tilde{R}_{ijk0}-\tilde{\Gamma}_{0i}^{p}\tilde{R}_{pjk0}-\tilde{\Gamma}%
_{0j}^{p}\tilde{R}_{ipk0}-\tilde{\Gamma}_{0k}^{p}\tilde{R}_{ijp0}%
-\tilde{\Gamma}_{00}^{p}\tilde{R}_{ijkp}\\
&  =t\left(  \frac{\partial}{\partial t}P_{ijk}+R_{i}^{p}P_{pjk}+R_{j}%
^{p}P_{ipk}+R_{k}^{p}P_{ijp}+\frac{1}{2}R_{ijkp}\nabla^{p}R\right) \\
&  +2P_{ijk}%
\end{align*}
and
\[
\tilde{\nabla}_{q}\tilde{R}_{ijk0}=\nabla_{q}P_{ijk}-\tilde{\Gamma}_{q0}%
^{p}\tilde{R}_{ijkp}=\nabla_{q}P_{ijk}+R_{q}^{p}R_{ijkp}+\frac{1}{2t}%
R_{ijkq}.
\]
Noticing that $\nabla^{q}R_{ijkq}=P_{ijk}$ by the second Bianchi identity, we
write the diffusion term in the form%
\begin{align*}
\tilde{\Delta}\tilde{R}_{ijk0}  &  =tg^{pq}\tilde{\nabla}_{p}\tilde{\nabla
}_{q}\tilde{R}_{ijk0}\\
&  =t\nabla^{q}\left(  \nabla_{q}P_{ijk}+R_{q}^{p}R_{ijkp}+\frac{1}%
{2t}R_{ijkq}\right) \\
&  +g^{pq}\left(  tR_{p}^{m}+\frac{1}{2}\delta_{p}^{m}\right)  \left(
\nabla_{q}R_{ijkm}\right) \\
&  =t\left(  \Delta P_{ijk}+\frac{1}{2}R_{ijkp}\nabla^{p}R+2R_{q}^{p}%
\nabla^{q}R_{ijkp}\right)  +P_{ijk}%
\end{align*}
and cancel terms to obtain%
\begin{align*}
\tilde{\nabla}_{0}\tilde{R}_{ijk0}  &  =\tilde{\Delta}\tilde{R}_{ijk0}%
+t\left[  \left(  \frac{\partial}{\partial t}-\Delta\right)  P_{ijk}+R_{p}%
^{q}\vee_{q}^{p}P_{ijk}-2R_{p}^{q}\nabla_{q}R_{ijk}^{p}\right] \\
&  +P_{ijk}.
\end{align*}
Thus the special case $\tilde{\nabla}_{0}\tilde{R}_{iik0}=\digamma_{ijk0}$ of
equation (\ref{st-riem-40-b-evol}) holds if and only if%
\[
D_{t}P_{ijk}=\Delta P_{ijk}+2R_{p}^{q}\nabla_{q}R_{ijk}^{p}+2g^{pq}\left(
R_{pij}^{r}P_{qrk}-R_{pjk}^{r}P_{qir}+R_{pik}^{r}P_{qjr}\right)  ,
\]
hence if and only if (\ref{P-evolve}) holds.

Finally, the equivalence of (\ref{R-evolve}) and the case $i,j,k,\ell\geq1$ of
(\ref{st-riem-40-b-evol}) is clear when we observe that%
\[
\digamma_{ijk\ell}=\Delta R_{ijk\ell}+2\left(  B_{ijk\ell}-B_{jik\ell
}-B_{jki\ell}+B_{ikj\ell}\right)  +\frac{1}{t}R_{ijk\ell}%
\]
and%
\[
\tilde{\nabla}_{0}\tilde{R}_{ijk\ell}=t\frac{\partial}{\partial t}\left(
\frac{1}{t}R_{ijk\ell}\right)  +R_{p}^{q}\vee_{q}^{p}R_{ijk\ell}+\frac{2}%
{t}R_{ijk\ell}=D_{t}R_{ijk\ell}+\frac{1}{t}R_{ijk\ell}.
\]

\end{proof}

\chapter{Generalized space-time connections\label{Generalize}}

In this section, we derive new matrix LYH inequalities for the Ricci flow by
generalizing the definition of the space-time connection in
\S \ \ref{Glossary}.

So let $\left(  \mathcal{M}^{n},\bar{g}\left(  \bar{t}\right)  \right)  $ be a
solution of the Ricci flow rescaled by a cosmological constant $\mu$:%
\[
\frac{\partial}{\partial\bar{t}}\bar{g}=-2\left(  \overline{\operatorname*{Rc}%
}+\mu\bar{g}\right)  .
\]
Consider the family of symmetric connections $\tilde{\nabla}$ defined on
space-time $(\widetilde{\mathcal{M}},\tilde{g})$ by%
\begin{align}
\tilde{\Gamma}_{ij}^{k}  &  \doteqdot\bar{\Gamma}_{ij}^{k} \tag{GC1}%
\label{GC1}\\
\tilde{\Gamma}_{i0}^{k}  &  \doteqdot-\left(  \bar{R}_{i}^{k}+\mu\delta
_{i}^{k}+A_{i}^{k}\right) \tag{GC2}\label{GC2}\\
\tilde{\Gamma}_{00}^{k}  &  \doteqdot-\left(  \frac{1}{2}\bar{\nabla}^{k}%
\bar{R}+B^{k}\right) \tag{GC3}\label{GC3}\\
\tilde{\Gamma}_{00}^{0}  &  \doteqdot-\left(  \mu+C\right)  \tag{GC4}%
\label{GC4}\\
\tilde{\Gamma}_{ij}^{0}  &  \doteqdot\tilde{\Gamma}_{i0}^{0}=0, \tag{GC5}%
\label{GC5}%
\end{align}
for $i,j,k\geq1$, where $A$ is a tensor of type $\left(  1,1\right)  $, $B$ is
a vector field, and $C$ is a scalar function. We saw in \S \ref{Glossary} that
the space-time connection $\tilde{\nabla}$ has a number of useful and
interesting properties when $A=B=C=0$. Our goal here is to investigate what
conditions on $A$, $B$, and $C$ are necessary and sufficient for
$\tilde{\nabla}$ to retain certain desirable characteristics. In particular,
we determine which connections of this form are both compatible with the
space-time metric and satisfy the Ricci flow for degenerate metrics. Such
space-time connections are worth studying, because their curvatures satisfy
parabolic evolution equations and thus furnish Li--Yau--Hamilton quadratics
for the Ricci flow.

Define a $\left(  2,0\right)  $-tensor $\bar{A}$ by%
\[
\bar{A}_{ij}\doteqdot A_{i}^{p}\bar{g}_{pj}.
\]
Our first observation is that $\tilde{\nabla}$ is both torsion-free and
compatible with $\tilde{g}$ exactly when $\bar{A}$ is a $2$-form:

\begin{lemma}
\label{parallel}The metric $\tilde{g}$ is parallel with respect to the
symmetric connection $\tilde{\nabla}$,%
\[
\tilde{\nabla}_{i}\tilde{g}^{jk}=0,
\]
if and only if (\ref{GC1})-(\ref{GC5}) hold, where $\bar{A}$ is a $2$-form,%
\begin{equation}
A_{p}^{j}\bar{g}^{pk}+A_{p}^{k}\bar{g}^{jp}=0, \label{A-2-form}%
\end{equation}
and there are no restrictions on either $B$ or $C$.
\end{lemma}

\begin{proof}
For $i,j,k\geq1$, the equation
\[
0=\tilde{\nabla}_{i}\tilde{g}^{jk}=\partial_{i}\tilde{g}^{jk}+\tilde{\Gamma
}_{ip}^{j}\tilde{g}^{pk}+\tilde{\Gamma}_{ip}^{k}\tilde{g}^{jp}%
\]
is equivalent to%
\[
\tilde{\Gamma}_{ij}^{k}=\bar{\Gamma}_{ij}^{k},
\]
since $\bar{\nabla}$ is the unique torsion-free connection compatible with
$\bar{g}$; this is (\ref{GC1}). Assuming (\ref{GC2}), the equation%
\[
0=\tilde{\nabla}_{0}\tilde{g}^{jk}=\partial_{0}\tilde{g}^{jk}+\tilde{\Gamma
}_{0p}^{j}\tilde{g}^{pk}+\tilde{\Gamma}_{0p}^{k}\tilde{g}^{jp}%
\]
is valid for $j,k\geq1$ if and only if%
\[
A_{p}^{j}\bar{g}^{pk}+A_{p}^{k}\bar{g}^{jp}=0.
\]
This says that when we lower an index, $\bar{A}_{ij}=A_{i}^{p}\bar{g}_{pj}$ is
a $2$-form. The equation
\[
0=\tilde{\nabla}_{i}\tilde{g}^{0k}=\partial_{i}\tilde{g}^{0k}+\tilde{\Gamma
}_{ip}^{0}\tilde{g}^{pk}+\tilde{\Gamma}_{ip}^{k}\tilde{g}^{0p}%
\]
is valid for $i\geq0$ and $k\geq1$ if and only if%
\[
\tilde{\Gamma}_{ip}^{0}=0
\]
holds for all $i\geq0$ and $p\geq1$; this is (\ref{GC5}). The identity
$\tilde{\nabla}_{i}\tilde{g}^{00}$ is satisfied automatically for all $i\geq0$.
\end{proof}

Hence by lowering indices, we may regard $\bar{A}$ as a $2$-form, $\bar{B}%
_{i}\doteqdot B^{p}\bar{g}_{pi}$ as a $1$-form, and $C$ as a $0$-form.

\section{The Riemann curvature tensor}

The space-time Riemann curvature tensor is defined by (\ref{Rm-def}); in
components, one has%
\begin{equation}
\tilde{R}_{ijk}^{\ell}=\partial_{i}\tilde{\Gamma}_{jk}^{\ell}-\partial
_{j}\tilde{\Gamma}_{ik}^{\ell}+\tilde{\Gamma}_{jk}^{m}\tilde{\Gamma}%
_{im}^{\ell}-\tilde{\Gamma}_{ik}^{m}\tilde{\Gamma}_{jm}^{\ell}.
\label{Rm-components}%
\end{equation}
By definition, we have the asymmetry $\tilde{R}_{ijk}^{\ell}=-\tilde{R}%
_{jik}^{\ell}$. The remaining formulas are as follows:

\begin{proposition}
If $i,j,k,\ell\geq1$ and $a,b,c\geq0$, then $\widetilde{\operatorname*{Rm}}$
satisfies:%
\begin{align}
\tilde{R}_{ijk}^{\ell}  &  =\bar{R}_{ijk}^{\ell}\tag{GR1}\label{GR1}\\
\tilde{R}_{ij0}^{\ell}  &  =\bar{\nabla}_{j}\bar{R}_{i}^{\ell}-\bar{\nabla
}_{i}\bar{R}_{j}^{\ell}+\bar{\nabla}_{j}A_{i}^{\ell}-\bar{\nabla}_{i}%
A_{j}^{\ell}\tag{GR2a}\label{GR2a}\\
\tilde{R}_{0jk}^{\ell}  &  =\bar{\nabla}^{\ell}\bar{R}_{jk}-\bar{\nabla}%
_{k}\bar{R}_{j}^{\ell}+\bar{\nabla}_{j}A_{k}^{\ell}\tag{GR2b}\label{GR2b}\\
\tilde{R}_{0j0}^{\ell}  &  =-\frac{\partial}{\partial\bar{t}}\bar{R}_{j}%
^{\ell}+\left(  \mu-C\right)  \bar{R}_{j}^{\ell}+\frac{1}{2}\bar{\nabla}%
_{j}\bar{\nabla}^{\ell}\bar{R}+\bar{R}_{j}^{m}\bar{R}_{m}^{\ell}%
\tag{GR3}\label{GR3}\\
&  -\frac{\partial}{\partial\bar{t}}A_{j}^{\ell}+A_{j}^{m}A_{m}^{\ell}+\left(
\mu-C\right)  A_{j}^{\ell}+\bar{R}_{j}^{m}A_{m}^{\ell}+A_{j}^{m}\bar{R}%
_{m}^{\ell}\nonumber\\
&  +\bar{\nabla}_{j}B^{\ell}-\mu C\delta_{j}^{\ell}.\nonumber\\
\tilde{R}_{abc}^{0}  &  =0. \tag{GR4}\label{GR4}%
\end{align}

\end{proposition}

\begin{proof}
Identities (\ref{GR1}) and (\ref{GR4}) follow easily from (\ref{Rm-components}%
). To derive (\ref{GR2a}), we use (\ref{GC2}) to compute%
\begin{align*}
\tilde{R}_{ij0}^{\ell}  &  =\partial_{i}\tilde{\Gamma}_{j0}^{\ell}%
+\tilde{\Gamma}_{im}^{\ell}\tilde{\Gamma}_{j0}^{m}-\partial_{j}\tilde{\Gamma
}_{i0}^{\ell}-\tilde{\Gamma}_{jm}^{\ell}\tilde{\Gamma}_{i0}^{m}\\
&  =\bar{\nabla}_{j}\left(  \bar{R}_{i}^{\ell}+A_{i}^{\ell}\right)
-\bar{\nabla}_{i}\left(  \bar{R}_{j}^{\ell}+A_{j}^{\ell}\right)  .
\end{align*}
To derive (\ref{GR2b}), we recall that%
\[
\frac{\partial}{\partial\bar{t}}\bar{\Gamma}_{jk}^{\ell}=-\bar{\nabla}_{j}%
\bar{R}_{k}^{\ell}-\bar{\nabla}_{k}\bar{R}_{j}^{\ell}+\bar{\nabla}^{\ell}%
\bar{R}_{jk}%
\]
and calculate%
\begin{align*}
\tilde{R}_{0jk}^{\ell}  &  =\partial_{0}\tilde{\Gamma}_{jk}^{\ell}-\left(
\partial_{j}\tilde{\Gamma}_{0k}^{\ell}-\tilde{\Gamma}_{jk}^{m}\tilde{\Gamma
}_{0m}^{\ell}+\tilde{\Gamma}_{jm}^{\ell}\tilde{\Gamma}_{0k}^{m}\right) \\
&  =\frac{\partial}{\partial\bar{t}}\bar{\Gamma}_{jk}^{\ell}+\bar{\nabla}%
_{j}\left(  \bar{R}_{k}^{\ell}+A_{k}^{\ell}\right) \\
&  =\bar{\nabla}^{\ell}\bar{R}_{jk}-\bar{\nabla}_{k}\bar{R}_{j}^{\ell}%
+\bar{\nabla}_{j}A_{k}^{\ell}.
\end{align*}
Finally, to derive (\ref{GR3}), we use (\ref{GC3}) and (\ref{GC4}) to compute%
\begin{align*}
\tilde{R}_{0j0}^{\ell}  &  =\partial_{0}\tilde{\Gamma}_{j0}^{\ell}-\left(
\partial_{j}\tilde{\Gamma}_{00}^{\ell}+\tilde{\Gamma}_{jm}^{\ell}\tilde
{\Gamma}_{00}^{m}\right)  +\tilde{\Gamma}_{j0}^{m}\tilde{\Gamma}_{0m}^{\ell}\\
&  =-\frac{\partial}{\partial\bar{t}}\left(  \bar{R}_{j}^{\ell}+A_{j}^{\ell
}\right)  +\bar{\nabla}_{j}\left(  \frac{1}{2}\bar{\nabla}^{\ell}\bar
{R}+B^{\ell}\right) \\
&  -\left(  \mu+C\right)  \left(  \bar{R}_{j}^{\ell}+\mu\delta_{j}^{\ell
}+A_{j}^{\ell}\right) \\
&  +\left(  \bar{R}_{j}^{m}+\mu\delta_{j}^{m}+A_{j}^{m}\right)  \left(
\bar{R}_{m}^{\ell}+\mu\delta_{m}^{\ell}+A_{m}^{\ell}\right) \\
&  =-\frac{\partial}{\partial\bar{t}}\bar{R}_{j}^{\ell}+\left(  \mu-C\right)
\bar{R}_{j}^{\ell}+\frac{1}{2}\bar{\nabla}_{j}\bar{\nabla}^{\ell}\bar{R}%
+\bar{R}_{j}^{m}\bar{R}_{m}^{\ell}\\
&  -\frac{\partial}{\partial\bar{t}}A_{j}^{\ell}+A_{j}^{m}A_{m}^{\ell}+\left(
\mu-C\right)  A_{j}^{\ell}+\bar{R}_{j}^{m}A_{m}^{\ell}+A_{j}^{m}\bar{R}%
_{m}^{\ell}\\
&  +\bar{\nabla}_{j}B^{\ell}-\mu C\delta_{j}^{\ell}.
\end{align*}

\end{proof}

\begin{corollary}
\label{gen-Ricci}If $i,j\geq1$, then $\widetilde{\operatorname*{Rc}}$
satisfies:%
\begin{align*}
\tilde{R}_{ij}  &  =\bar{R}_{ij}\\
\tilde{R}_{0k}  &  =\frac{1}{2}\bar{\nabla}_{k}\bar{R}-\left(  \bar{\delta
}\bar{A}\right)  _{k}\\
\tilde{R}_{00}  &  =\frac{1}{2}\frac{\partial}{\partial\bar{t}}\bar
{R}+C\left(  \bar{R}+n\mu\right)  +\left|  \bar{A}\right|  _{\bar{g}}^{2}%
+\bar{\delta}\bar{B},
\end{align*}
where
\[
\left(  \bar{\delta}\bar{A}\right)  _{k}\doteqdot-\bar{\nabla}^{p}\bar{A}%
_{pk}=\bar{\nabla}_{p}A_{k}^{p}%
\]
and
\[
\bar{\delta}\bar{B}\doteqdot-\bar{\nabla}^{p}\bar{B}_{p}=-\bar{\nabla}%
_{p}B^{p}.
\]

\end{corollary}

\begin{proof}
The first two equations are easy. For the third, we substitute the formula%
\[
\frac{1}{2}\frac{\partial}{\partial\bar{t}}\bar{R}=\frac{1}{2}\bar{\Delta}%
\bar{R}+\bar{R}_{pq}\bar{R}^{pq}+\mu\bar{R}%
\]
into the calculation%
\begin{align*}
\tilde{R}_{00}  &  =-\tilde{R}_{0j0}^{j}\\
&  =\frac{\partial}{\partial\bar{t}}\bar{R}+\left(  C-\mu\right)  \bar
{R}-\frac{1}{2}\bar{\Delta}\bar{R}-\bar{R}_{j}^{m}\bar{R}_{m}^{m}-A_{j}%
^{m}A_{m}^{j}-\bar{\nabla}_{j}B^{j}+n\mu C
\end{align*}
and cancel terms.
\end{proof}

\section[Ricci flow for degenerate metrics]{Solutions of the Ricci flow for
degenerate metrics}

The goal of this section is to determine necessary and sufficient conditions
on $A$, $B$, and $C$ for $\left(  \tilde{g},\tilde{\nabla}\right)  $ to
satisfy the rescaled Ricci flow for degenerate metrics. (Recall Definition
\ref{degenerateRF}.) Our results here are most easily stated if we introduce
the $1$-form%
\[
\bar{E}\doteqdot\bar{B}+2\bar{\delta}\bar{A}.
\]
We shall see that a particularly nice set of equations is obtained when
$\bar{A}$ and $\bar{E}$ are closed initially. In this case, there is always a
solution $\left(  \tilde{g},\tilde{\nabla}\right)  $ satisfying Definition
\ref{degenerateRF} for as long as $\bar{g}\left(  \bar{t}\right)  $ exists.

\begin{proposition}
\label{Christoffel-evolve}Suppose $C=\mu$. Then $\left(  \tilde{g}%
,\tilde{\nabla}\right)  $ satisfy the Ricci flow with cosmological term $\mu$,
namely%
\begin{equation}
\frac{\partial}{\partial\bar{t}}\tilde{\Gamma}_{ij}^{k}=-\tilde{\nabla}%
_{i}\tilde{R}_{j}^{k}-\tilde{\nabla}_{j}\tilde{R}_{i}^{k}+\tilde{\nabla}%
^{k}\tilde{R}_{ij}, \label{degen-RF-eqn}%
\end{equation}
if and only if the $2$-form $\bar{A}=\bar{A}_{ij}\,dx^{i}\otimes dx^{j}$
satisfies%
\begin{equation}
\frac{\partial}{\partial\bar{t}}\bar{A}=-d\bar{\delta}\bar{A}-2\mu\bar{A}
\label{Abar-eqn}%
\end{equation}
and the $1$-form $\bar{E}=\bar{E}_{i}dx^{i}$ satisfies%
\begin{equation}
\frac{\partial}{\partial\bar{t}}\bar{E}=-d\bar{\delta}\bar{E}-2\mu\bar
{E}-d\left|  \bar{A}\right|  _{\bar{g}}^{2}. \label{Ebar-eqn}%
\end{equation}
If $d\bar{A}=0$ initially, then $d\bar{A}\equiv0$ for as long as a solution
exists; and if $d\bar{E}=0$ initially, then $d\bar{E}\equiv0$ for as long as a
solution exists. So if $\bar{A}$ and $\bar{E}$ are closed initially,
(\ref{degen-RF-eqn}) is valid if and only if $\bar{A}$ and $\bar{E}$ evolve by%
\begin{equation}
\frac{\partial}{\partial\bar{t}}\bar{A}=\bar{\Delta}_{d}\bar{A}-2\mu\bar{A}
\label{A-closed-evolve}%
\end{equation}
and%
\begin{equation}
\frac{\partial}{\partial\bar{t}}\bar{E}=\bar{\Delta}_{d}\bar{E}-2\mu\bar
{E}-d\left|  \bar{A}\right|  _{\bar{g}}^{2} \label{E-closed-evolve}%
\end{equation}
respectively, where $-\bar{\Delta}_{d}\doteqdot d\bar{\delta}+\bar{\delta}d$
is the Hodge--de Rham Laplacian.
\end{proposition}

\begin{remark}
When $\bar{A}$ and $\bar{E}$ are closed initially, (\ref{A-closed-evolve}) and
(\ref{E-closed-evolve}) are both parabolic equations whose solutions exist as
long as the solution of the Ricci flow with cosmological constant $\mu$ exists.
\end{remark}

\begin{remark}
The choice $C=\mu$ is useful to obtain good evolution equation if either $A$
or $B$ is nonzero. But if $A$ and $B$ are both identically zero, taking $C=0$
as in \S \ref{Glossary} generally yields better results.
\end{remark}

\begin{remark}
\label{scaling-remark}If $\left(  \bar{A},\bar{E}\right)  $ is a pair of
initially-closed forms satisfying equations (\ref{A-closed-evolve}) and
(\ref{E-closed-evolve}), then the pair $\left(  \lambda\bar{A},\lambda^{2}%
\bar{E}\right)  $ is also, for any $\lambda\in\mathbb{R}$.
\end{remark}

\begin{proof}
[Proof of Proposition \ref{Christoffel-evolve}]Let $\digamma_{ij}^{k}$ denote
the RHS of (\ref{degen-RF-eqn}). If $i,j,k\geq1$, then formula
(\ref{degen-RF-eqn}) reduces to the standard evolution equation for
$\bar{\Gamma}_{ij}^{k}$. It is easily checked that both sides of
(\ref{degen-RF-eqn}) vanish if $k=0$, provided that $\mu$ and $C$ are
constant. If $j=0$ but $i,k\geq1$, then%
\begin{align*}
\digamma_{i0}^{k}  &  \doteqdot-\tilde{\nabla}_{i}\tilde{R}_{0}^{k}%
-\tilde{\nabla}_{0}\tilde{R}_{i}^{k}+\tilde{\nabla}^{k}\tilde{R}_{i0}\\
&  =-\bar{\nabla}_{i}\left[  \frac{1}{2}\bar{\nabla}^{k}\bar{R}-\left(
\bar{\delta}\bar{A}\right)  ^{k}\right]  -\left(  \bar{R}_{i}^{p}+\mu
\delta_{i}^{p}+A_{i}^{p}\right)  \bar{R}_{p}^{k}\\
&  -\frac{\partial}{\partial\bar{t}}\bar{R}_{i}^{k}-\left(  \bar{R}_{i}%
^{p}+\mu\delta_{i}^{p}+A_{i}^{p}\right)  \bar{R}_{p}^{k}+\left(  \bar{R}%
_{p}^{k}+\mu\delta_{p}^{k}+A_{p}^{k}\right)  \bar{R}_{i}^{p}\\
&  +\bar{\nabla}^{k}\left[  \frac{1}{2}\bar{\nabla}_{i}\bar{R}-\left(
\bar{\delta}\bar{A}\right)  _{i}\right]  +\bar{g}^{k\ell}\left(  \bar{R}%
_{\ell}^{p}+\mu\delta_{\ell}^{p}+A_{\ell}^{p}\right)  \bar{R}_{ip}\\
&  =-\frac{\partial}{\partial\bar{t}}\bar{R}_{i}^{k}-2A_{i}^{p}\bar{R}_{p}%
^{k}-\bar{\nabla}_{i}\bar{\nabla}^{p}A_{p}^{k}-\bar{\nabla}^{k}\bar{\nabla
}_{p}A_{i}^{p}.
\end{align*}
Since $\tilde{\Gamma}_{i0}^{k}=-\left(  \bar{R}_{i}^{k}+\mu\delta_{i}%
^{k}+A_{i}^{k}\right)  $, it follows that (\ref{degen-RF-eqn}) holds for $j=0$
and $i,k\geq1$ if and only if%
\[
\frac{\partial}{\partial\bar{t}}A_{i}^{k}=\bar{\nabla}_{i}\bar{\nabla}%
^{p}A_{p}^{k}+\bar{\nabla}^{k}\bar{\nabla}_{p}A_{i}^{p}+2A_{i}^{p}\bar{R}%
_{p}^{k},
\]
hence if and only if%
\[
\frac{\partial}{\partial\bar{t}}\bar{A}_{ij}=\frac{\partial}{\partial\bar{t}%
}\left(  A_{i}^{k}\bar{g}_{kj}\right)  =-\left(  d\bar{\delta}\bar{A}\right)
_{ij}-2\mu\bar{A}_{ij}.
\]
If $i=j=0$ but $k\geq1$, we recall that%
\[
\bar{\nabla}^{k}\left(  \frac{\partial}{\partial\bar{t}}\bar{R}\right)
=\frac{\partial}{\partial\bar{t}}\left(  \bar{\nabla}^{k}\bar{R}\right)
-2\left(  \bar{R}_{\ell}^{k}\bar{\nabla}^{\ell}\bar{R}+\mu\bar{\nabla}^{k}%
\bar{R}\right)
\]
and compute%
\begin{align*}
\digamma_{00}^{k}  &  \doteqdot-\tilde{\nabla}_{0}\tilde{R}_{0}^{k}%
-\tilde{\nabla}_{0}\tilde{R}_{0}^{k}+\tilde{\nabla}^{k}\tilde{R}_{00}\\
&  =-2\left[  \frac{\partial}{\partial\bar{t}}\left(  \frac{1}{2}\bar{\nabla
}^{k}\bar{R}-\left(  \bar{\delta}\bar{A}\right)  ^{k}\right)  -\tilde{\Gamma
}_{00}^{p}\tilde{R}_{p}^{k}+\tilde{\Gamma}_{0p}^{k}\tilde{R}_{0}^{p}\right] \\
&  +\bar{\nabla}^{k}\left(  \frac{1}{2}\frac{\partial}{\partial\bar{t}}\bar
{R}+C\bar{R}+\left|  \bar{A}\right|  _{\bar{g}}^{2}+\bar{\delta}\bar
{B}\right)  -2\bar{g}^{k\ell}\tilde{\Gamma}_{\ell0}^{p}\tilde{R}_{p0}\\
&  =-\frac{1}{2}\frac{\partial}{\partial\bar{t}}\left(  \bar{\nabla}^{k}%
\bar{R}\right)  -\bar{\nabla}^{k}\bar{\nabla}_{p}B^{p}-2B^{p}\bar{R}_{p}^{k}\\
&  +2\frac{\partial}{\partial\bar{t}}\left(  \bar{g}^{k\ell}\bar{\nabla}%
_{p}A_{\ell}^{p}\right)  +\bar{\nabla}^{k}\left|  \bar{A}\right|  _{\bar{g}%
}^{2}+2\left(  \mu-C\right)  \left(  \bar{\nabla}^{p}A_{p}^{k}\right)
-4\bar{R}^{kp}\bar{\nabla}_{q}A_{p}^{q}.
\end{align*}
Since $\tilde{\Gamma}_{00}^{k}=-\frac{1}{2}\bar{\nabla}^{k}\bar{R}-B^{k}$, it
follows that (\ref{degen-RF-eqn}) holds for $i=j=0$ and $k\geq1$ if and only
if%
\begin{align*}
\frac{\partial}{\partial\bar{t}}\left(  \bar{B}_{j}+2\bar{\nabla}_{p}A_{j}%
^{p}\right)   &  =\frac{\partial}{\partial\bar{t}}\left[  \bar{g}_{jk}\left(
B^{k}+2\bar{g}^{k\ell}\bar{\nabla}_{p}A_{\ell}^{p}\right)  \right] \\
&  =\bar{\nabla}_{j}\bar{\nabla}_{p}B^{p}-2\mu\bar{B}_{j}-\bar{\nabla}%
_{j}\left|  \bar{A}\right|  _{\bar{g}}^{2}+2\left(  \mu+C\right)  \bar{\nabla
}^{p}\bar{A}_{pj}.
\end{align*}
When $C=\mu$, this equation is the same as%
\begin{align*}
\frac{\partial}{\partial\bar{t}}\left(  \bar{B}+2\bar{\delta}\bar{A}\right)
&  =-d\bar{\delta}\bar{B}-2\mu\left(  \bar{B}+2\bar{\delta}\bar{A}\right)
-d\left|  \bar{A}\right|  _{\bar{g}}^{2}\\
&  =-d\bar{\delta}\left(  \bar{B}+2\bar{\delta}\bar{A}\right)  -2\mu\left(
\bar{B}+2\bar{\delta}\bar{A}\right)  -d\left|  \bar{A}\right|  _{\bar{g}}^{2},
\end{align*}
because $\bar{\delta}^{2}=0$.

To complete the proof, it suffices to note that%
\[
\frac{\partial}{\partial\bar{t}}\left(  d\bar{A}\right)  =d\left(
\frac{\partial}{\partial\bar{t}}\bar{A}\right)  =-2\mu\left(  d\bar{A}\right)
\]
and%
\[
\frac{\partial}{\partial\bar{t}}\left(  d\bar{E}\right)  =d\left(
\frac{\partial}{\partial\bar{t}}\bar{E}\right)  =-2\mu\left(  d\bar{E}\right)
,
\]
because the exterior derivative is independent of the metric and satisfies
$d^{2}=0$.
\end{proof}

In analogy with Remark \ref{cov-Ricci-sym}, we make the following observation:

\begin{remark}
If $\bar{A}$ evolves according to equation (\ref{Abar-eqn}), one has the
symmetry%
\[
\tilde{\nabla}_{i}\tilde{R}_{j0}-\tilde{\nabla}_{j}\tilde{R}_{i0}%
=\tilde{\nabla}_{0}\bar{A}_{ij}.
\]

\end{remark}

\begin{proof}
Because%
\begin{align*}
\tilde{\nabla}_{i}\tilde{R}_{j0}  &  =\bar{\nabla}_{i}\left(  \frac{1}{2}%
\bar{\nabla}_{j}\bar{R}+\bar{\nabla}^{p}\bar{A}_{pj}\right)  -\tilde{\Gamma
}_{i0}^{p}\bar{R}_{jp}\\
&  =\frac{1}{2}\bar{\nabla}_{i}\bar{\nabla}_{j}\bar{R}+\bar{\nabla}_{i}%
\bar{\nabla}^{p}\bar{A}_{pj}+\bar{R}_{i}^{p}\bar{R}_{pj}+\mu\bar{R}_{ij}%
+\bar{A}_{i}^{p}\bar{R}_{pj},
\end{align*}
we observe that when (\ref{Abar-eqn}) holds, we have%
\begin{align*}
\tilde{\nabla}_{0}\bar{A}_{ij}  &  =\frac{\partial}{\partial\bar{t}}\bar
{A}_{ij}-\tilde{\Gamma}_{0i}^{p}\bar{A}_{pj}-\tilde{\Gamma}_{0j}^{p}\bar
{A}_{ip}\\
&  =\bar{\nabla}_{i}\bar{\nabla}^{p}\bar{A}_{pj}-\bar{\nabla}_{j}\bar{\nabla
}^{p}\bar{A}_{pi}+A_{i}^{p}\bar{R}_{pj}-A_{j}^{p}\bar{R}_{pi}\\
&  =\tilde{\nabla}_{i}\tilde{R}_{j0}-\tilde{\nabla}_{j}\tilde{R}_{i0}.
\end{align*}

\end{proof}

\section[New quadratics]{New Li--Yau--Hamilton quadratics}

We now wish to regard $\widetilde{\operatorname*{Rm}}$ as the bilinear form
defined on $\Lambda^{2}T\widetilde{\mathcal{M}}$ by (\ref{Bilinear1}) and
(\ref{Bilinear2}). To be useful as a LYH\ quadratic, it is desirable that a
bilinear form be symmetric and positive. Fortunately, symmetry of
$\widetilde{\operatorname*{Rm}}$ is compatible with the other properties we
wish $\tilde{\nabla}$ to possess. In particular, we have the following:

\begin{lemma}
The bilinear form $\widetilde{\operatorname*{Rm}}$ has the symmetry%
\[
\tilde{R}_{ij0\ell}=\tilde{R}_{0\ell ij}%
\]
for all $i,j,\ell\geq1$ if and only if $\bar{A}$ is a closed $2$-form.
Moreover, $\widetilde{\operatorname*{Rm}}$ has the symmetry%
\[
\tilde{R}_{0j0\ell}=\tilde{R}_{0\ell0j}%
\]
for all $j,\ell\geq1$ if $\bar{A}$ evolves by (\ref{Abar-eqn}), $C=\mu$, and
$\bar{E}$ is a closed $1$-form.
\end{lemma}

\begin{proof}
By (\ref{GR2a}) and (\ref{GR2b}), we have%
\[
\tilde{R}_{ij0\ell}-\tilde{R}_{0\ell ij}=\bar{g}_{\ell p}\tilde{R}_{ij0}%
^{p}-\bar{g}_{jp}\tilde{R}_{0\ell i}^{p}=\bar{\nabla}_{j}\bar{A}_{i\ell}%
-\bar{\nabla}_{i}\bar{A}_{j\ell}+\bar{\nabla}_{\ell}\bar{A}_{ij}=\left(
d\bar{A}\right)  _{ji\ell}.
\]
Next we observe that%
\[
\bar{g}_{jp}\frac{\partial}{\partial\bar{t}}\bar{R}_{\ell}^{p}-\bar{g}_{\ell
p}\frac{\partial}{\partial\bar{t}}\bar{R}_{j}^{p}=0
\]
and%
\begin{equation}
\bar{g}_{jp}\frac{\partial}{\partial\bar{t}}A_{\ell}^{p}-\bar{g}_{\ell p}%
\frac{\partial}{\partial\bar{t}}A_{j}^{p}=2\frac{\partial}{\partial\bar{t}%
}\bar{A}_{\ell j}+2\bar{R}_{jp}A_{\ell}^{p}-2\bar{R}_{\ell p}A_{j}^{p}%
+4\mu\bar{A}_{\ell j}. \label{lower-commute}%
\end{equation}
Hence by (\ref{GR3}),%
\begin{align*}
\tilde{R}_{0j0\ell}-\tilde{R}_{0\ell0j}  &  =\bar{g}_{\ell p}\tilde{R}%
_{0j0}^{p}-\bar{g}_{jp}\tilde{R}_{0\ell0}^{p}\\
&  =\bar{g}_{jp}\left(  \frac{\partial}{\partial\bar{t}}\bar{R}_{\ell}%
^{p}+\frac{\partial}{\partial\bar{t}}A_{\ell}^{p}\right)  -\bar{g}_{\ell
p}\left(  \frac{\partial}{\partial\bar{t}}\bar{R}_{j}^{p}+\frac{\partial
}{\partial\bar{t}}A_{j}^{p}\right) \\
&  +2\left(  \mu-C\right)  \bar{A}_{j\ell}+2A_{j}^{k}\bar{R}_{k\ell}-2A_{\ell
}^{k}\bar{R}_{kj}+\left(  \bar{\nabla}_{j}\bar{B}_{\ell}-\bar{\nabla}_{\ell
}\bar{B}_{j}\right) \\
&  =2\frac{\partial}{\partial\bar{t}}\bar{A}_{\ell j}+2\left(  \mu+C\right)
\bar{A}_{\ell j}-\left(  d\bar{B}\right)  _{\ell j}.
\end{align*}
If $\frac{\partial}{\partial\bar{t}}\bar{A}=-d\bar{\delta}\bar{A}-2\mu\bar{A}%
$, this becomes%
\[
\tilde{R}_{0j0\ell}-\tilde{R}_{0\ell0j}=2\left(  C-\mu\right)  \bar{A}_{\ell
j}-d\left(  \bar{B}+2\bar{\delta}\bar{A}\right)  _{\ell j}.
\]

\end{proof}

It is also fortunate that the maximum principle applies to the curvature
$\widetilde{\operatorname*{Rm}}$ of a generalized connection. To see this, it
will be convenient to introduce a $\left(  1,1\right)  $-tensor $\tilde{A}$
defined for $i,j\geq1$ by%
\begin{align*}
\tilde{A}_{i}^{j}  &  =A_{i}^{j}\\
\tilde{A}_{0}^{j}  &  =\left(  B+\bar{\delta}\bar{A}\right)  ^{j}\\
\tilde{A}_{i}^{0}  &  =0\\
\tilde{A}_{0}^{0}  &  =\mu.
\end{align*}

\begin{proposition}
\label{gen-Rm-evolve}Let $\bar{A}$ and $\bar{E}$ be closed initially and
evolve by (\ref{A-closed-evolve}) and (\ref{E-closed-evolve}), respectively.
Let $C=\mu$ be constant. Then $\widetilde{\operatorname*{Rm}}$ is a symmetric
bilinear form which evolves by%
\begin{equation}
\tilde{\nabla}_{0}\widetilde{\operatorname*{Rm}}=\tilde{\Delta}\widetilde
{\operatorname*{Rm}}+\widetilde{\operatorname*{Rm}}^{2}+\widetilde
{\operatorname*{Rm}}^{\#}+2\mu\widetilde{\operatorname*{Rm}}+\tilde{A}%
\vee\widetilde{\operatorname*{Rm}}. \label{gen-Rm-evolve-eq1}%
\end{equation}

\end{proposition}

\begin{proof}
Because $\tilde{\nabla}$ is symmetric, one computes directly from the
definition that%
\begin{align*}
\frac{\partial}{\partial\bar{t}}\tilde{R}_{ijk}^{\ell}  &  =\partial
_{i}\left(  \frac{\partial}{\partial t}\tilde{\Gamma}_{jk}^{\ell}\right)
-\partial_{j}\left(  \frac{\partial}{\partial t}\tilde{\Gamma}_{ik}^{\ell
}\right) \\
&  +\left(  \frac{\partial}{\partial t}\tilde{\Gamma}_{jk}^{m}\right)
\tilde{\Gamma}_{im}^{\ell}+\tilde{\Gamma}_{jk}^{m}\left(  \frac{\partial
}{\partial t}\tilde{\Gamma}_{im}^{\ell}\right)  -\left(  \frac{\partial
}{\partial t}\tilde{\Gamma}_{ik}^{m}\right)  \tilde{\Gamma}_{jm}^{\ell}%
-\tilde{\Gamma}_{ik}^{m}\left(  \frac{\partial}{\partial t}\tilde{\Gamma}%
_{jm}^{\ell}\right) \\
&  =\partial_{i}\left(  \frac{\partial}{\partial t}\tilde{\Gamma}_{jk}^{\ell
}\right)  -\tilde{\Gamma}_{ij}^{m}\left(  \frac{\partial}{\partial t}%
\tilde{\Gamma}_{mk}^{\ell}\right)  -\tilde{\Gamma}_{ik}^{m}\left(
\frac{\partial}{\partial t}\tilde{\Gamma}_{jm}^{\ell}\right)  +\tilde{\Gamma
}_{im}^{\ell}\left(  \frac{\partial}{\partial t}\tilde{\Gamma}_{jk}^{m}\right)
\\
&  -\partial_{j}\left(  \frac{\partial}{\partial t}\tilde{\Gamma}_{ik}^{\ell
}\right)  +\tilde{\Gamma}_{ji}^{m}\left(  \frac{\partial}{\partial t}%
\tilde{\Gamma}_{mk}^{\ell}\right)  +\tilde{\Gamma}_{jk}^{m}\left(
\frac{\partial}{\partial t}\tilde{\Gamma}_{im}^{\ell}\right)  -\tilde{\Gamma
}_{jm}^{\ell}\left(  \frac{\partial}{\partial t}\tilde{\Gamma}_{ik}^{m}\right)
\\
&  =\tilde{\nabla}_{i}\left(  \frac{\partial}{\partial t}\tilde{\Gamma}%
_{jk}^{\ell}\right)  -\tilde{\nabla}_{j}\left(  \frac{\partial}{\partial
t}\tilde{\Gamma}_{ik}^{\ell}\right)  .
\end{align*}
Hence by Proposition \ref{Christoffel-evolve} and the Ricci identities,%
\begin{align*}
\frac{\partial}{\partial\bar{t}}\tilde{R}_{ijk}^{\ell}  &  =\bar{\nabla}%
_{j}\left(  \nabla_{k}\tilde{R}_{i}^{\ell}-\tilde{\nabla}^{\ell}\tilde{R}%
_{ik}\right)  -\tilde{\nabla}_{i}\left(  \tilde{\nabla}_{k}\tilde{R}_{j}%
^{\ell}-\tilde{\nabla}^{\ell}\tilde{R}_{jk}\right) \\
&  +\tilde{R}_{jim}^{\ell}\tilde{R}_{k}^{m}-\tilde{R}_{jik}^{m}\tilde{R}%
_{m}^{\ell}.
\end{align*}
On the other hand, the second Bianchi identity implies that%
\begin{align*}
\tilde{\Delta}\tilde{R}_{ijk}^{\ell}  &  \doteqdot\tilde{g}^{pq}\tilde{\nabla
}_{p}\tilde{\nabla}_{q}\tilde{R}_{ijk}^{\ell}=-\tilde{g}^{pq}\tilde{\nabla
}_{p}\left(  \tilde{\nabla}_{i}\tilde{R}_{jqk}^{\ell}+\tilde{\nabla}_{j}%
\tilde{R}_{qik}^{\ell}\right) \\
&  =\tilde{\nabla}_{i}\tilde{\nabla}^{\ell}\tilde{R}_{jk}-\tilde{\nabla}%
_{i}\tilde{\nabla}_{k}\tilde{R}_{j}^{\ell}-\tilde{\nabla}_{j}\tilde{\nabla
}^{\ell}\tilde{R}_{ik}\\
&  +\tilde{\nabla}_{j}\tilde{\nabla}_{k}\tilde{R}_{i}^{\ell}-\tilde{R}_{i}%
^{m}\tilde{R}_{jmk}^{\ell}-\tilde{R}_{j}^{m}\tilde{R}_{mik}^{\ell}\\
&  +\tilde{g}^{pq}\left[
\begin{array}
[c]{c}%
\tilde{R}_{pij}^{m}\tilde{R}_{mqk}^{\ell}+\tilde{R}_{pik}^{m}\tilde{R}%
_{jqm}^{\ell}-\tilde{R}_{pim}^{\ell}\tilde{R}_{jqk}^{m}\\
\mathstrut\\
+\tilde{R}_{pji}^{m}\tilde{R}_{qmk}^{\ell}+\tilde{R}_{pjk}^{m}\tilde{R}%
_{qim}^{\ell}-\tilde{R}_{pjm}^{\ell}\tilde{R}_{qik}^{m}%
\end{array}
\right]  ,
\end{align*}
while straightforward calculations reveal that%
\[
\tilde{R}_{ijk\ell}^{2}=\tilde{g}^{pq}\tilde{g}^{rs}\tilde{R}_{ijpr}\tilde
{R}_{sqk\ell}=-\tilde{g}^{pq}\left(  \tilde{R}_{pij}^{m}\tilde{R}_{mqk\ell
}+\tilde{R}_{pji}^{m}\tilde{R}_{qmk\ell}\right)
\]
and%
\begin{align*}
\tilde{R}_{ijk\ell}^{\#}  &  =\tilde{R}_{abcd}\tilde{R}_{pqrs}\left(
\delta_{i}^{a}\delta_{j}^{q}\tilde{g}^{bp}-\delta_{i}^{p}\delta_{j}^{b}%
\tilde{g}^{aq}\right)  \left(  \delta_{\ell}^{c}\delta_{k}^{s}\tilde{g}%
^{dr}-\delta_{\ell}^{r}\delta_{k}^{d}\tilde{g}^{cs}\right) \\
&  =-\tilde{g}^{pq}\left(  \tilde{R}_{qim\ell}\tilde{R}_{pjk}^{m}+\tilde
{R}_{pik}^{m}\tilde{R}_{jqm\ell}+\tilde{R}_{jpm\ell}\tilde{R}_{qik}^{m}%
+\tilde{R}_{jqk}^{m}\tilde{R}_{ipm\ell}\right)  .
\end{align*}
Because equation (\ref{gen-Rm-evolve-eq1}) is readily verified for $k=\ell=0$,
we may assume without loss of generality that $\ell\geq1$. Then we can combine
the identities above to obtain%
\begin{align*}
\frac{\partial}{\partial\bar{t}}\tilde{R}_{ijk\ell}  &  =\bar{g}_{\ell m}%
\frac{\partial}{\partial\bar{t}}\tilde{R}_{ijk}^{m}-2\bar{R}_{\ell m}\tilde
{R}_{ijk}^{m}-2\mu\tilde{R}_{ijk\ell}\\
&  =\tilde{\Delta}\tilde{R}_{ijk\ell}+\tilde{R}_{ijk\ell}^{2}+\tilde
{R}_{ijk\ell}^{\#}-2\mu\tilde{R}_{ijk\ell}-\tilde{R}_{p}^{q}\vee_{q}^{p}%
\tilde{R}_{ijk\ell},
\end{align*}
because $\bar{g}_{lm}\tilde{\Delta}\tilde{R}_{ijk}^{m}=\tilde{\Delta}\tilde
{R}_{ijk\ell}$ and $\tilde{R}_{\ell}^{m}\tilde{R}_{ijkm}=\bar{R}_{\ell
m}\tilde{R}_{ijk}^{m}$. Now if we regard $\tilde{\Gamma}_{0}$ as a
(globally-defined) space-time $\left(  1,1\right)  $-tensor $\left(
\tilde{\Gamma}_{0}\right)  _{p}^{q}=\tilde{\Gamma}_{0p}^{q}$, we may write%
\[
\tilde{\nabla}_{0}\tilde{R}_{ijk\ell}=\frac{\partial}{\partial\bar{t}}%
\tilde{R}_{ijk\ell}-\left(  \tilde{\Gamma}_{0}\vee\widetilde
{\operatorname*{Rm}}\right)  _{ijk\ell},
\]
yielding%
\begin{equation}
\tilde{\nabla}_{0}\widetilde{\operatorname*{Rm}}=\tilde{\Delta}\widetilde
{\operatorname*{Rm}}+\widetilde{\operatorname*{Rm}}^{2}+\widetilde
{\operatorname*{Rm}}^{\#}-2\mu\widetilde{\operatorname*{Rm}}-\left(
\tilde{\Gamma}_{0}+\widetilde{\operatorname*{Rc}}\right)  \vee\widetilde
{\operatorname*{Rm}}, \label{gen-Rm-evolve-eq2}%
\end{equation}
where $\widetilde{\operatorname*{Rc}}$ here denotes the $\left(  1,1\right)
$-tensor $\tilde{R}_{p}^{q}\doteqdot\tilde{R}_{pm}\tilde{g}^{mq}$. We claim
equation (\ref{gen-Rm-evolve-eq2}) is equivalent to equation
(\ref{gen-Rm-evolve-eq1}). Indeed, it follows from (\ref{GC2})--(\ref{GC4})
and Corollary \ref{gen-Ricci} that for $p,q\geq1$ one has%
\begin{align*}
\left(  \tilde{\Gamma}_{0}+\widetilde{\operatorname*{Rc}}\right)  _{p}^{q}  &
=-\mu\delta_{p}^{q}-A_{p}^{q}\\
\left(  \tilde{\Gamma}_{0}+\widetilde{\operatorname*{Rc}}\right)  _{0}^{q}  &
=-B^{q}-\left(  \bar{\delta}\bar{A}\right)  ^{q}\\
\left(  \tilde{\Gamma}_{0}+\widetilde{\operatorname*{Rc}}\right)  _{p}^{0}  &
=0\\
\left(  \tilde{\Gamma}_{0}+\widetilde{\operatorname*{Rc}}\right)  _{0}^{0}  &
=-2\mu.
\end{align*}
So one need only check that if $N$ denotes the number of space-like components
of $\tilde{R}_{ijk\ell}$, the RHS of (\ref{gen-Rm-evolve-eq2}) contains
$-2+N+2\left(  4-N\right)  =6-N$ terms of the form $\mu\tilde{R}_{ijk\ell}$,
while the RHS of (\ref{gen-Rm-evolve-eq1}) contains $2+\left(  4-N\right)
=6-N$ such terms.
\end{proof}

\bigskip

Armed with the tools to show that $\widetilde{\operatorname*{Rm}}$ remains
non-negative, we are now ready to construct Li--Yau--Hamilton quadratics.

\begin{condition}
Assume in the remainder of this paper that $g\left(  t\right)  $ is a solution
of the Ricci flow on $\mathcal{M}$ for $t\in\lbrack0,\Omega)$, and that
$\bar{g}\left(  \bar{t}\right)  \equiv e^{-\bar{t}}g\left(  e^{\bar{t}%
}\right)  $ is the associated solution of (\ref{cosmo}), the Ricci flow with
cosmological constant $\mu=1/2$, for $\bar{t}\in\left(  -\infty,\ln
\Omega\right)  $. Assume further that the generalized connection on
$(\widetilde{\mathcal{M}},\tilde{g})$ defined by (\ref{GC1})--(\ref{GC5}) has
$C=\mu=1/2$.
\end{condition}

Given any $2$-form $U$ and $1$-form $W$ on $\mathcal{M}$, we define $\tilde
{X}\doteqdot U\oplus\frac{1}{2}W$, so that for $i,j,\geq1$,%
\begin{align*}
\tilde{X}^{ij}  &  =U^{ij}\\
\tilde{X}^{0j}  &  =-\tilde{X}^{j0}=\frac{1}{2}W^{j}.
\end{align*}
If there exist $\bar{A}$ and $\bar{B}$ such that $\widetilde
{\operatorname*{Rm}}$ is symmetric, we define the forms%
\begin{align*}
A_{ik}  &  \doteqdot A_{i}^{j}g_{jk}=t\bar{A}_{ik}\\
B_{k}  &  \doteqdot B^{j}g_{jk}=t\bar{B}_{k}\\
E_{k}  &  \doteqdot\left(  B^{j}+2\left(  \bar{\delta}\bar{A}\right)
^{j}\right)  g_{jk}=B_{k}+2t\left(  \delta A\right)  _{k},
\end{align*}
and make the following observations:

\begin{remark}
\label{tbar-t}Let $\bar{\Delta}_{d}$ and $\Delta_{d}$ denote the Hodge--de
Rham Laplacians of $\left(  \bar{g},\bar{\nabla}\right)  $ and $\left(
g,\nabla\right)  $, respectively. Then $\bar{A}$ is closed and evolves by%
\[
\frac{\partial}{\partial\bar{t}}\bar{A}=\bar{\Delta}_{d}\bar{A}-2\mu\bar{A}
\]
if and only if $A$ is closed and evolves by%
\[
\frac{\partial}{\partial t}A=\Delta_{d}A+\frac{1-2\mu}{t}A.
\]
Moreover, $\bar{E}\equiv\bar{B}+2\bar{\delta}\bar{A}$ is closed and evolves by%
\[
\frac{\partial}{\partial\bar{t}}\bar{E}=\bar{\Delta}_{d}\bar{E}-2\mu\bar
{E}-d\left|  \bar{A}\right|  _{\bar{g}}^{2}%
\]
if and only if $E\equiv B+2t\delta A$ is closed and evolves by%
\[
\frac{\partial}{\partial t}E=\Delta_{d}E+\frac{1-2\mu}{t}E-d\left|  A\right|
_{g}^{2}.
\]
If $\mu=1/2$ and $A$ is closed (exact) initially, then $A$ remains closed
(exact). Likewise, if $\mu=1/2$ and $E$ is closed (exact) initially, then $E$
remains closed (exact).
\end{remark}

\begin{proof}
The evolution equations are straightforward calculations. Together with
Proposition \ref{Christoffel-evolve}, they imply that $A$ and $E$ remain
closed if they are initially. To prove the assertions about exactness, suppose
that $A\left(  0\right)  =d\alpha_{0}$ and $E\left(  0\right)  =d\varepsilon
_{0}$. Let $\alpha\left(  t\right)  $ and $\varepsilon\left(  t\right)  $ be
solutions of%
\[
\frac{\partial}{\partial t}\alpha=\Delta_{d}\alpha,\quad\quad\quad\quad
\alpha\left(  0\right)  =\alpha_{0}%
\]
and%
\[
\frac{\partial}{\partial t}\varepsilon=\Delta_{d}\varepsilon-\left|  A\right|
_{g}^{2},\quad\quad\quad\quad\varepsilon\left(  0\right)  =\varepsilon_{0}%
\]
respectively. Then%
\[
\frac{\partial}{\partial t}\left(  d\alpha\right)  =d\left(  \frac{\partial
}{\partial t}\alpha\right)  =-d\left(  d\delta+\delta d\right)  \alpha
=-\left(  d\delta+\delta d\right)  \left(  d\alpha\right)  =\Delta_{d}\left(
d\alpha\right)  ,
\]
and similarly%
\[
\frac{\partial}{\partial t}\left(  d\varepsilon\right)  =d\left(
\frac{\partial}{\partial t}\varepsilon\right)  =d\left(  \Delta_{d}%
\varepsilon-\left|  A\right|  _{g}^{2}\right)  =\Delta_{d}\left(
d\varepsilon\right)  -d\left|  A\right|  _{g}^{2}.
\]
By uniqueness of solutions to parabolic equations, we have $A\left(  t\right)
\equiv d\alpha\left(  t\right)  $ and $E\left(  t\right)  \equiv
d\varepsilon\left(  t\right)  $ for as long as $g\left(  t\right)  $ exists.
\end{proof}

\begin{theorem}
\label{full}Let $\left(  \mathcal{M}^{n},g\left(  t\right)  \right)  $ be a
solution of the Ricci flow on a closed manifold and a time interval
$[0,\Omega)$. Let $A_{0}$ be a $2$-form which is closed at $t=0$ and let
$E_{0}$ be a $1$-form which is closed at $t=0$. Then there is a solution
$A\left(  t\right)  $ of%
\[
\frac{\partial}{\partial t}A=\Delta_{d}A,\quad\quad\quad\quad A\left(
0\right)  =A_{0}%
\]
and a solution $E\left(  t\right)  $ of%
\[
\frac{\partial}{\partial t}E=\Delta_{d}E-d\left\vert A\right\vert _{g}%
^{2},\quad\quad\quad\quad E\left(  0\right)  =E_{0}%
\]
which exist for all $t\in\lbrack0,\Omega)$. Suppose that%
\[
0\leq\operatorname*{Rm}\left(  U,U\right)  +\frac{1}{4}\left\vert W\right\vert
^{2}+\left\vert A\left(  W\right)  \right\vert ^{2}-2\left\langle \nabla
_{W}A,U\right\rangle -\left\langle \nabla_{W}E,W\right\rangle
\]
at $t=0$ for any $2$-form $U$ and $1$-form $W$ on $\mathcal{M}$, and let
$\widetilde{\operatorname*{Rm}}$ be the curvature of the generalized
connection on $(\widetilde{\mathcal{M}},\tilde{g})$ with $B\equiv E-2t\delta
A$ and $C=1/2$. Then for all $t\in\lbrack0,\Omega)$, one has the estimate%
\begin{align*}
0  &  \leq e^{\bar{t}}\widetilde{\operatorname*{Rm}}\left(  \tilde{X}%
,\tilde{X}\right) \\
&  =R_{ijk\ell}U^{ij}U^{\ell k}+2\left[  t\left(  \nabla_{\ell}R_{jk}%
-\nabla_{k}R_{j\ell}\right)  +\nabla_{j}A_{k\ell}\right]  W^{j}U^{\ell k}\\
&  +\left[
\begin{array}
[c]{c}%
t^{2}\left(  \Delta R_{j\ell}-\frac{1}{2}\nabla_{j}\nabla_{\ell}R+2R_{jpq\ell
}R^{pq}-R_{j}^{p}R_{p\ell}\right) \\
+t\left(  2\mu R_{j\ell}+A_{j}^{p}R_{p\ell}+A_{\ell}^{p}R_{pj}-2\nabla
_{j}\nabla^{k}A_{k\ell}\right) \\
+\mu^{2}g_{j\ell}-A_{j}^{p}A_{p\ell}-\nabla_{j}E_{\ell}%
\end{array}
\right]  W^{j}W^{\ell}.
\end{align*}

\end{theorem}

\begin{proof}
By Remark \ref{tbar-t}, there are closed solutions $\bar{A}$ and $\bar{E}$ of
(\ref{A-closed-evolve}) and (\ref{E-closed-evolve}), respectively, existing
for $-\infty<\bar{t}<\log\Omega$. Set $Q\doteqdot e^{\bar{t}}\widetilde
{\operatorname*{Rm}}$. Then by Proposition \ref{gen-Rm-evolve},%
\begin{align}
\frac{\partial}{\partial t}Q  &  =\frac{d\bar{t}}{dt}\frac{\partial}%
{\partial\bar{t}}\left(  e^{\bar{t}}\widetilde{\operatorname*{Rm}}\right)
\nonumber\\
&  =\frac{\partial}{\partial\bar{t}}\widetilde{\operatorname*{Rm}}%
+\widetilde{\operatorname*{Rm}}\nonumber\\
&  =\tilde{\Delta}\widetilde{\operatorname*{Rm}}+\widetilde{\operatorname*{Rm}%
}^{2}+\widetilde{\operatorname*{Rm}}^{\#}+2\widetilde{\operatorname*{Rm}%
}+\left(  \tilde{\Gamma}_{0}+\tilde{A}\right)  \vee\widetilde
{\operatorname*{Rm}}. \label{Q-evolve}%
\end{align}
Since $Q\geq0$ at $t=0$, it will suffice to apply the space-time maximum
principle (Proposition \ref{STMP}) to $Q$ for $0\leq t<\Omega$. Notice that
for $j,k\geq1$,%
\begin{align*}
\left(  \tilde{\Gamma}_{0}+\tilde{A}\right)  _{j}^{k}  &  =-tR_{j}^{k}%
-\frac{1}{2}\delta_{j}^{k}\\
\left(  \tilde{\Gamma}_{0}+\tilde{A}\right)  _{0}^{k}  &  =-\frac{t^{2}}%
{2}\nabla^{k}R+t\left(  \delta A\right)  ^{k}\\
\left(  \tilde{\Gamma}_{0}+\tilde{A}\right)  _{j}^{0}  &  =0\\
\left(  \tilde{\Gamma}_{0}+\tilde{A}\right)  _{0}^{0}  &  =-\frac{1}{2},
\end{align*}
and define a generalized symmetric connection $\hat{\nabla}$ on $\mathcal{M}%
\times\lbrack0,\Omega)$ for $i,j,k\geq1$ by%
\begin{align*}
\hat{\Gamma}_{ij}^{k}  &  \doteqdot\Gamma_{ij}^{k}\\
\hat{\Gamma}_{i0}^{k}  &  \doteqdot-R_{i}^{k}\\
\hat{\Gamma}_{00}^{k}  &  \doteqdot-\left(  \frac{1}{2}\nabla^{k}%
R+F^{k}\right) \\
\hat{\Gamma}_{00}^{0}  &  =\hat{\Gamma}_{i0}^{0}=\hat{\Gamma}_{ij}^{0}=0,
\end{align*}
where the $\Gamma_{ij}^{k}$ are determined by the Levi-Civita connection of
$g$, and
\[
F^{k}\doteqdot\frac{1}{2}\left(  t-1\right)  \nabla^{k}R-\left(  \delta
A\right)  ^{k}.
\]
Let $\hat{g}$ be the space-time metric on $\mathcal{M}\times\lbrack0,\Omega)$
induced by $g$: namely, $\hat{g}^{ij}=g^{ij}$ and $\hat{g}^{i0}=\hat{g}%
^{00}=0$ for $i,j\geq1$. By Lemma \ref{parallel}, $\hat{g}$ is parallel with
respect to $\hat{\nabla}$. One verifies by straightforward calculation that%
\[
\tilde{\Delta}\widetilde{\operatorname*{Rm}}\doteqdot\tilde{g}^{pq}%
\tilde{\nabla}_{p}\tilde{\nabla}_{q}\widetilde{\operatorname*{Rm}}=\hat
{g}^{pq}\hat{\nabla}_{p}\hat{\nabla}_{q}Q\doteqdot\hat{\Delta}Q
\]
and
\[
\tilde{R}_{ijk\ell}^{2}=\tilde{g}^{ad}\tilde{g}^{bc}\tilde{R}_{ijab}\tilde
{R}_{cdk\ell}=g^{ad}g^{bc}Q_{ijab}Q_{cdk\ell}\doteqdot Q_{ijk\ell}^{2},
\]
and%
\[
\tilde{R}_{ijk\ell}^{\#}=\tilde{R}_{abcd}\tilde{R}_{pqrs}C_{ij}^{ab,pq}C_{\ell
k}^{cd,rs}=Q_{abcd}Q_{pqrs}c_{ij}^{ab,pq}c_{\ell k}^{cd,rs}\doteqdot
Q_{ijk\ell}^{\#},
\]
where $c_{ij}^{ab,cd}\doteqdot\delta_{i}^{a}\delta_{j}^{d}g^{bc}-\delta
_{i}^{c}\delta_{j}^{b}g^{ad}=e^{-\bar{t}}C_{ij}^{ab,cd}$. Thus equation
(\ref{Q-evolve}) can be written in the form%
\[
\hat{\nabla}_{0}Q=\hat{\Delta}Q+Q^{2}+Q^{\#},
\]
whence the theorem follows by applying Proposition \ref{STMP} to $Q$ for
$0\leq t<\Omega.$
\end{proof}

If $\widetilde{\operatorname*{Rm}}$ is symmetric, we define another symmetric
bilinear form $\tilde{H}$ on $\Lambda^{2}T\widetilde{\mathcal{M}}$ by
stipulating that $\tilde{H}$ obey the symmetries
\[
\tilde{H}_{abcd}=-\tilde{H}_{bacd}=-\tilde{H}_{abdc}=\tilde{H}_{cdab}%
\]
for all $a,b,c,d\geq0$ and satisfy%
\begin{align*}
\tilde{H}_{ijk\ell}  &  \doteqdot\bar{R}_{ijk\ell}\\
\tilde{H}_{0jk\ell}  &  \doteqdot\bar{\nabla}_{j}\bar{A}_{k\ell}\\
\tilde{H}_{0j0\ell}  &  \doteqdot A_{j}^{p}\bar{A}_{p\ell}+\bar{\nabla}%
_{j}\bar{E}_{\ell}%
\end{align*}
for all $i,j,k,\ell\geq1$. Then we define and expand%
\begin{align*}
\Psi\left(  A,E,U,W\right)   &  \doteqdot e^{\bar{t}}\tilde{H}\left(
\tilde{X},\tilde{X}\right) \\
&  =R_{ijk\ell}U^{ij}U^{\ell k}+2W^{j}\nabla_{j}A_{k\ell}U^{\ell k}\\
&  -\left(  A_{j}^{p}A_{p\ell}+\nabla_{j}E_{\ell}\right)  W^{j}W^{\ell}.
\end{align*}

\begin{theorem}
\label{main}Let $\left(  \mathcal{M}^{n},g\left(  t\right)  \right)  $ be a
solution of the Ricci flow on a closed manifold and a time interval
$[0,\Omega)$. Let $A_{0}$ be a $2$-form which is closed at $t=0$, and let
$E_{0}$ be a $1$-form which is closed at $t=0$. Then there is a solution
$A\left(  t\right)  $ of%
\[
\frac{\partial}{\partial t}A=\Delta_{d}A,\quad\quad\quad\quad A\left(
0\right)  =A_{0}%
\]
and a solution $E\left(  t\right)  $ of%
\[
\frac{\partial}{\partial t}E=\Delta_{d}E-d\left|  A\right|  _{g}^{2}%
,\quad\quad\quad\quad E\left(  0\right)  =E_{0}%
\]
which exist for all $t\in\lbrack0,\Omega)$. Suppose that%
\begin{align*}
0  &  \leq\left.  \Psi\left(  A,E,U,W\right)  \right|  _{t=0}\\
&  =\operatorname*{Rm}\left(  U,U\right)  -2\left\langle \nabla_{W}%
A,U\right\rangle +\left|  A\left(  W\right)  \right|  ^{2}-\left\langle
\nabla_{W}E,W\right\rangle
\end{align*}
for any $2$-form $U$ and $1$-form $W$ on $\mathcal{M}$. Then $\Psi\left(
A,E,U,W\right)  $ remains non-negative for all $t\in\lbrack0,\Omega)$.
\end{theorem}

\begin{proof}
Suppose $\left(  \mathcal{M},g\left(  t\right)  \right)  $ exists for
$t\in\lbrack0,\Omega)$, and define
\[
\Phi\left(  A,E,U,W\right)  \doteqdot e^{\bar{t}}\widetilde{\operatorname*{Rm}%
}\left(  \tilde{X},\tilde{X}\right)  ,
\]
where $\widetilde{\operatorname*{Rm}}$ is the curvature of the generalized
connection on $(\widetilde{\mathcal{M}},\tilde{g})$ with $B\equiv E-2t\delta
A$ and $C=1/2$. By hypothesis,
\[
0\leq\Psi\left(  A,E,U,W\right)  \leq\Phi\left(  A,E,U,W\right)
\]
at $t=0$ for all $U$ and $W$. So $\Phi\left(  A,E,U,W\right)  \geq0$ for all
$t\in\lbrack0,\Omega)$ and all $U,W$ by Theorem \ref{full}. Now since
$\Psi\left(  \lambda A,\lambda^{2}E,U,W\right)  =\Psi\left(  A,E,U,\lambda
W\right)  $, it follows from Remark \ref{scaling-remark} that $\Phi\left(
\lambda A,\lambda^{2}E,U,W\right)  \geq0$ for all $t\in\lbrack0,\Omega)$, all
$U,W$, and all $\lambda>0$. In particular, at each fixed $t\in\lbrack
0,\Omega)$, we have $0\leq\Phi\left(  \lambda A,\lambda^{2}E,U,\lambda
^{-1}W\right)  $ for all $U,W$ and $\lambda>0$, and hence%
\[
0\leq\lim_{\lambda\rightarrow\infty}\Phi\left(  \lambda A,\lambda
^{2}E,U,\lambda^{-1}W\right)  =\Psi\left(  A,E,U,W\right)  .
\]

\end{proof}

\begin{corollary}
\label{trace-Harnack}Let the hypotheses of Theorem \ref{main} hold. Then for
any $1$-form $V$, we have%
\begin{align*}
0\leq\psi\left(  A,E,V\right)   &  \doteqdot\operatorname*{Rc}\left(
V,V\right)  -2\left(  \delta A\right)  \left(  V\right)  +\left|  A\right|
^{2}+\delta E\\
&  =\operatorname*{Rc}\left(  V,V\right)  -2\left(  \delta A\right)  \left(
V\right)  +\left|  A\right|  ^{2}+\delta B
\end{align*}
for as long as the solution $\left(  \mathcal{M},g\left(  t\right)  \right)  $ exists.
\end{corollary}

\begin{proof}
For an orthonormal frame $\left\{  e_{k}\right\}  $, take $U_{ij}=\frac{1}%
{2}\left(  V_{i}W_{j}-V_{j}W_{i}\right)  $ and trace over $W\in\left\{
e_{1},\dots,e_{n}\right\}  $.
\end{proof}

\chapter{Examples\label{Apply}}

We shall now develop some examples, which may be regarded as further
corollaries of Theorem \ref{main}. Our intent is to explore the utility of the
new LYH quadratics by comparing their implications with known results in a few
special cases. Although our principle examples come from K\"{a}hler geometry,
we emphasize that the main results of this paper are more general, and in no
way require a K\"{a}hler structure.

\section{K\"{a}hler examples}

By definition, a Riemannian manifold $\left(  \mathcal{M}^{n},g\right)  $ is
K\"{a}hler if there is an almost-complex structure $J:T\mathcal{M}\rightarrow
T\mathcal{M}$ such that $g$ is Hermitian $g\left(  JX,JX\right)  =g\left(
X,X\right)  $ and $J$ is parallel $\nabla_{X}\left(  JY\right)  =J\left(
\nabla_{X}Y\right)  $. The latter condition immediately implies the symmetry
$R\left(  X,Y\right)  JZ=J\left(  R\left(  X,Y\right)  Z\right)  $ of the
Riemannian curvature.

There is a one-to-one correspondence between $J$-invariant symmetric
$2$-tensors and $J$-invariant $2$-forms on a K\"{a}hler manifold. Indeed,
given a $J$-invariant symmetric $2$-tensor $\phi\left(  JX,JY\right)
=\phi\left(  X,Y\right)  =\phi\left(  Y,X\right)  $, let $F\left(  X,Y\right)
\doteqdot\phi\left(  JX,Y\right)  $. Then $F$ is a $2$-form $F\left(
Y,X\right)  =-F\left(  X,Y\right)  $ which is $J$-invariant $F\left(
JX,JY\right)  =F\left(  X,Y\right)  $. Conversely, given a $J$-invariant
$2$-form $H\left(  JX,JY\right)  =H\left(  X,Y\right)  =-H\left(  Y,X\right)
$, let $\eta\left(  X,Y\right)  \doteqdot-H\left(  JX,Y\right)  $. Then $\eta$
is symmetric $\eta\left(  Y,X\right)  =\eta\left(  X,Y\right)  $ and
$J$-invariant $\eta\left(  JX,JY\right)  =\eta\left(  X,Y\right)  $. The
\emph{K\"{a}hler form }$\omega$ is defined to be the $2$-form induced by a
Hermitian metric $g$, and the \emph{Ricci form }$\rho$ is defined to be the
$2$-form induced by the Ricci tensor $\operatorname*{Rc}$ of $g$. These
definitions are justified by the following standard facts:

\begin{lemma}
Let $J$ be an almost-complex structure on $\left(  \mathcal{M}^{n},g\right)  $.

\begin{enumerate}
\item If $\left(  \mathcal{M},g\right)  $ is K\"{a}hler, then
$\operatorname*{Rc}$ is $J$-invariant.

\item If $\left(  \mathcal{M},g\right)  $ is K\"{a}hler, then $\rho$ is closed.

\item If $g$ is Hermitian, then $J$ is parallel$\ $if and only if $\omega$ is closed.
\end{enumerate}
\end{lemma}

It follows easily that the Ricci flow preserves a K\"{a}hler structure.
Indeed, (1) implies that $g$ remains Hermitian, whence (2) and (3) imply that
$J$ remains parallel, because%
\[
\frac{\partial}{\partial t}\left(  d\omega\right)  =d\left(  \frac{\partial
}{\partial t}\omega\right)  =d\left(  -2\rho\right)  =0.
\]
A key observation for the constructions in this section is the following:

\begin{lemma}
Let $\phi$ be a $J$-invariant symmetric $2$-tensor and let $F$ be the
corresponding $J$-invariant $2$-form. Then $\phi$ satisfies the
Lichnerowicz-Laplacian heat equation
\[
\frac{\partial}{\partial t}\phi=\Delta_{L}\phi,
\]
if and only if $F$ satisfies the Hodge-Laplacian heat equation
\[
\frac{\partial}{\partial t}F=\Delta_{d}F.
\]

\end{lemma}

\begin{proof}
If $\frac{\partial}{\partial t}\phi=\Delta_{L}\phi$, then direct computation
gives%
\begin{align*}
\frac{\partial}{\partial t}F_{ij}  &  =\frac{\partial}{\partial t}\left(
J_{i}^{k}\phi_{kj}\right)  =J_{i}^{k}\left(  \Delta\phi_{kj}+2R_{kpqj}%
\phi^{pq}-R_{j}^{\ell}\phi_{k\ell}-R_{k}^{\ell}\phi_{j\ell}\right) \\
&  =\Delta F_{ij}+2J_{i}^{k}R_{kpqj}\phi^{pq}-R_{j}^{\ell}F_{i\ell}-J_{i}%
^{k}R_{k}^{\ell}\phi_{j\ell}.
\end{align*}
The Hodge--de Rham Laplacian of $F$ is%
\begin{align*}
\Delta_{d}F_{ij}  &  =-\left(  d\delta F\right)  _{ij}-\left(  \delta
dF\right)  _{ij}\\
&  =\left(  \nabla_{i}\nabla^{k}F_{kj}-\nabla_{j}\nabla^{k}F_{ki}\right)
+\nabla^{k}\left(  \nabla_{k}F_{ij}-\nabla_{i}F_{kj}+\nabla_{j}F_{ki}\right)
\\
&  =\Delta F_{ij}-\left(  \nabla^{k}\nabla_{i}-\nabla_{i}\nabla^{k}\right)
F_{kj}+\left(  \nabla^{k}\nabla_{j}-\nabla_{j}\nabla^{k}\right)  F_{ki}\\
&  =\Delta F_{ij}+2g^{k\ell}R_{kij}^{m}F_{\ell m}-R_{i}^{k}F_{kj}-R_{j}%
^{k}F_{ik}.
\end{align*}
The identity
\begin{align*}
\left\langle R\left(  Z,W\right)  X,Y\right\rangle  &  =\left\langle R\left(
X,Y\right)  Z,W\right\rangle \\
&  =\left\langle R\left(  X,Y\right)  JZ,JW\right\rangle =\left\langle
R\left(  JZ,JW\right)  X,Y\right\rangle
\end{align*}
implies in particular that $R\left(  JZ,W\right)  X=-R\left(  Z,JW\right)  X$,
and hence%
\[
J_{i}^{k}R_{kpqj}\phi^{pq}=-R_{ikqj}J_{p}^{k}\phi^{pq}=g^{k\ell}%
R_{kijq}J_{\ell}^{p}\phi_{p}^{q}=g^{k\ell}R_{kij}^{q}F_{\ell q}.
\]
Thus%
\[
\frac{\partial}{\partial t}F_{ij}=\Delta_{d}F_{ij}-J_{i}^{k}R_{k}^{\ell}%
\phi_{j\ell}+R_{i}^{k}F_{kj}.
\]
But since the Ricci tensor is $J$-invariant, we have%
\[
J_{k}^{\ell}R_{i}^{k}=J_{i}^{k}R_{k}^{\ell},
\]
so that%
\[
-J_{i}^{k}R_{k}^{\ell}\phi_{j\ell}=-J_{k}^{\ell}R_{i}^{k}\phi_{j\ell}%
=-R_{i}^{k}F_{kj}.
\]
Hence $\frac{\partial}{\partial t}F=\Delta_{d}F$. The converse is proved similarly.
\end{proof}

\begin{proposition}
\label{full-kaehler}If $\left(  \mathcal{M}^{n},g\left(  t\right)  \right)  $
is a K\"{a}hler solution of the Ricci flow with non-negative curvature
operator on a closed manifold, the choices
\[
A=t\rho+\frac{1}{2}\omega
\]
and
\[
E=-\frac{t^{2}}{2}\,dR
\]
yield the estimate%
\begin{align*}
0  &  \leq\operatorname*{Rm}\left(  U,U\right)  -2\left\langle \nabla_{W}%
\rho,U\right\rangle +\frac{1}{4t^{2}}\left|  W\right|  ^{2}+\frac{1}%
{t}\operatorname*{Rc}\left(  W,W\right) \\
&  +\operatorname*{Rc}{}^{2}\left(  W,W\right)  +\frac{1}{2}\left(
\nabla\nabla R\right)  \left(  W,W\right)
\end{align*}
for all $t>0$ such that the solution exists.
\end{proposition}

\begin{proof}
The choice $A=t\rho+\frac{1}{2}\omega$ satisfies%
\[
\frac{\partial}{\partial t}A=\rho+t\frac{\partial}{\partial t}\rho+\frac{1}%
{2}\left(  -2\rho\right)  =t\Delta_{d}\rho=\Delta_{d}A
\]
by the preceding lemma, and
\[
\left|  A\right|  ^{2}=t^{2}\left|  \operatorname*{Rc}\right|  ^{2}%
+tR+\frac{n}{4}.
\]
If $E=df$ for some smooth function $f$, then $E$ will satisfy
\[
\frac{\partial}{\partial t}E=\Delta_{d}E-d\left|  A\right|  ^{2}%
\]
if and only if%
\[
d\left(  \frac{\partial}{\partial t}f\right)  =\frac{\partial}{\partial
t}E=d\left(  \Delta f-t^{2}\left|  \operatorname*{Rc}\right|  ^{2}-tR-\frac
{n}{4}\right)  .
\]
The choice $f=-t^{2}R/2$ satisfies%
\[
\frac{\partial}{\partial t}f=\Delta f-t^{2}\left|  \operatorname*{Rc}\right|
^{2}-tR.
\]
So to apply Theorem \ref{main}, we need only check that $\Psi\left(
A,E,U,W\right)  \geq0$ at $t=0$ for any $2$-form $U$ and $1$-form $W$. Noting
that $\operatorname*{Rc}{}^{2}\left(  JW,JW\right)  =\operatorname*{Rc}{}%
^{2}\left(  W,W\right)  $, we compute%
\begin{align*}
\Psi\left(  A,E,U,W\right)   &  =\operatorname*{Rm}\left(  U,U\right)
-2t\left\langle \nabla_{W}\rho,U\right\rangle +\frac{1}{4}\left|  W\right|
^{2}+t\operatorname*{Rc}\left(  W,W\right) \\
&  +t^{2}\operatorname*{Rc}{}^{2}\left(  W,W\right)  +\frac{t^{2}}{2}\left(
\nabla\nabla R\right)  \left(  W,W\right)  .
\end{align*}
This is clearly non-negative at $t=0$ whenever the curvature operator is.
Hence Theorem \ref{main} implies in particular that $0\leq\Psi\left(
A,E,U,\frac{1}{t}W\right)  $ for all $U$ and $W$ at any $t>0$ such that the
solution exists.
\end{proof}

To apply Corollary \ref{trace-Harnack}, we note that $\delta E=\delta
df=\frac{t^{2}}{2}\Delta R$. Because%
\[
\left(  \delta A\right)  _{i}=t\left(  \delta\rho\right)  _{i}+\frac{1}%
{2}\left(  \delta\omega\right)  _{i}=t\nabla^{j}\rho_{ij}=tJ_{i}^{k}\nabla
^{j}R_{kj}=\frac{t}{2}J_{i}^{k}\nabla_{k}R,
\]
we have
\[
\left(  \delta A\right)  \left(  V\right)  =\frac{t}{2}J_{i}^{k}\nabla
_{k}RV^{i}=\frac{t}{2}\left\langle \nabla R,JV\right\rangle .
\]
Then setting $V=tJX$, we compute%
\begin{align*}
\psi\left(  A,E,V\right)   &  =\operatorname*{Rc}\left(  V,V\right)
-t\left\langle \nabla R,JV\right\rangle \\
&  +\left(  t^{2}\left|  \operatorname*{Rc}\right|  ^{2}+tR+\frac{n}%
{4}\right)  +\frac{t^{2}}{2}\Delta R\\
&  =\frac{t^{2}}{2}\left[
\begin{array}
[c]{c}%
\Delta R+2\left|  \operatorname*{Rc}\right|  ^{2}+\frac{R}{t}\\
+2\left\langle \nabla R,X\right\rangle +2\operatorname*{Rc}\left(
JX,JX\right)
\end{array}
\right] \\
&  +\frac{1}{2}\left(  tR+\frac{n}{2}\right)  ,
\end{align*}
obtaining:

\begin{example}
\label{trace-kaehler}If $\left(  \mathcal{M},g\left(  t\right)  \right)  $ is
a K\"{a}hler solution of the Ricci flow with non-negative curvature operator
on a closed manifold, then for any vector field $X$ and all $t>0$ such that a
solution exists, we have the estimate%
\begin{equation}
0\leq\left[  \frac{\partial}{\partial t}R+\frac{R}{t}+2\left\langle \nabla
R,X\right\rangle +2\operatorname*{Rc}\left(  X,X\right)  \right]  +\left(
\frac{R}{t}+\frac{n}{2t^{2}}\right)  . \label{trace-kaehler-inequality}%
\end{equation}

\end{example}

The terms in square brackets compose Hamilton's trace quadratic. So this
special case of our general inequality is weaker than that estimate. To gauge
the potential usefulness of our inequality, however, we can make qualitative
comparisons. Taking $X=0$, we obtain%
\[
0\leq t^{2}\frac{\partial R}{\partial t}+2tR+\frac{n}{2}=\frac{\partial
}{\partial t}\left[  t\left(  tR+\frac{n}{2}\right)  \right]  .
\]
Hence at any $x\in\mathcal{M}$ and times $t_{2}\geq t_{1}$ with $t_{2}\neq0$,
we have%
\[
R\left(  x,t_{2}\right)  \geq\left(  \frac{t_{1}}{t_{2}}\right)  ^{2}R\left(
x,t_{1}\right)  +\frac{n}{2t_{2}}\left(  \frac{t_{1}}{t_{2}}-1\right)  .
\]
If we have an ancient solution, so that the interval of existence is
$t\in\left(  -\infty,\Omega\right)  $, we can translate in time $t\mapsto
t-\tau$ and take the limit as $\tau\rightarrow\infty$ to conclude that $R$ is
a pointwise nondecreasing function of $t$. Thus this particular example
(\ref{trace-kaehler-inequality}) of our general result is strong enough to
recover that important fact. (Compare \cite{HamEternal} and \cite{HamForm}.)

On the other hand, suppose $t_{1}>0$ and there is some constant $C>0$ such
that $R\left(  x,t_{1}\right)  \geq nC/t_{1}$. Then for all $t_{2}\in\left[
t_{1},Ct_{1}\right]  $ we have%
\[
R\left(  x,t_{2}\right)  \geq\frac{R\left(  x,t_{1}\right)  }{2C^{2}}+\left(
\frac{t_{1}}{t_{2}}\right)  ^{2}\frac{nC}{2t_{1}}+\frac{n}{2t_{2}}\left(
\frac{1}{C}-1\right)  \geq\frac{R\left(  x,t_{1}\right)  }{2C^{2}}.
\]
In all known applications of a trace LYH inequality, one has $t_{1}\geq c>0$.
Thus in this sense also, the special case (\ref{trace-kaehler-inequality}) of
Corollary \ref{trace-Harnack} is qualitatively equivalent to Hamilton's estimate.

\begin{remark}
\label{no-time}If $\left(  \mathcal{M}^{n},g\left(  t\right)  \right)  $ is a
K\"{a}hler solution of the Ricci flow on a closed manifold, we can choose
$A=\rho$ and $E=-dR/2$ (thereby dropping the explicit time dependency from
Example \ref{trace-kaehler}) and thus obtain the quadratic%
\begin{align*}
\Psi\left(  \rho,-\frac{1}{2}dR,U,W\right)   &  =\operatorname*{Rm}\left(
U,U\right)  -2\left\langle \nabla_{W}\rho,U\right\rangle \\
&  +\operatorname*{Rc}{}^{2}\left(  W,W\right)  +\frac{1}{2}\left(
\nabla\nabla R\right)  \left(  W,W\right)  .
\end{align*}
One would not, however, expect this to be positive for general initial data.
\end{remark}

\section{Other examples}

\begin{proposition}
\label{surface-system}Let $\left(  \mathcal{M}^{2},g\left(  t\right)  \right)
$ be a solution of the Ricci flow on a closed surface, and let $dS$ denote the
area element of $g$. Then for any pair $\left(  \phi,f\right)  $ solving the
system%
\begin{subequations}
\begin{align}
\frac{\partial}{\partial t}\phi &  =\Delta\phi+R\phi\label{surface-system-a}\\
\frac{\partial}{\partial t}f  &  =\Delta f+\phi^{2}, \label{surface-system-b}%
\end{align}
the choices $A=\phi\,dS$ and $E=-2df$ yield the trace quadratic%
\end{subequations}
\[
\psi\left(  A,E,X\right)  =R\left\vert X\right\vert ^{2}+2\left\langle
\nabla\phi,X\right\rangle +\frac{\partial}{\partial t}f
\]
and the matrix quadratic%
\begin{align*}
&  \Psi\left(  A,E,U,W\right) \\
&  =R\left\vert U\right\vert ^{2}-2\left\langle \nabla\phi,W\right\rangle
\left\langle \omega,U\right\rangle +\phi^{2}\left\vert W\right\vert
^{2}+2\left(  \nabla\nabla f\right)  \left(  W,W\right)  .
\end{align*}
If $\psi\geq0$ $\left(  \Psi\geq0\right)  $ at $t=0$, then $\psi\geq0$
$\left(  \Psi\geq0\right)  $ for as long as the solution exists.
\end{proposition}

\begin{example}
\label{surface-trace}On any closed surface $\left(  \mathcal{M}^{2},g\right)
$ of non-negative curvature evolving by the Ricci flow, the pair $\left(
\phi,f\right)  $ given by%
\begin{align*}
\phi &  =tR+1\\
f  &  =t\phi=t^{2}R+t
\end{align*}
yield the estimate%
\[
0\leq\frac{\partial}{\partial t}R+2\left\langle \nabla R,X\right\rangle
+R\left|  X\right|  ^{2}+\frac{2tR+1}{t^{2}}%
\]
for any vector field $X$ and all $t>0$ such that a solution exists.
\end{example}

It is remarkable that this special case of Proposition \ref{surface-system}
recovers the result in Example \ref{trace-kaehler}. Even though Proposition
\ref{surface-system} is \emph{a priori }more general than Proposition
\ref{full-kaehler} in the sense that it depends only on the dimension, the
construction in Proposition \ref{surface-system} is specific to surfaces, and
does not generalize readily to higher-dimension manifolds (whether K\"{a}hler
or not).

In analogy with Remark \ref{no-time}, one can drop the explicit time
dependency in Example \ref{surface-trace} and notice that the choices
$\phi=f=R$ solve system (\ref{surface-system-a})--(\ref{surface-system-b}), obtaining:

\begin{remark}
On any closed surface $\left(  \mathcal{M}^{2},g\left(  t\right)  \right)  $
evolving by the Ricci flow, one has the trace LYH quadratic
\[
\psi\left(  A,E,X\right)  =\frac{\partial}{\partial t}R+2\left\langle \nabla
R,X\right\rangle +R\left\vert X\right\vert ^{2}.
\]

\end{remark}

\begin{proof}
[Proof of Proposition \ref{surface-system}]Write the area element of $g$ as
$\left(  dS\right)  _{ij}=J_{i}^{k}g_{kj}$, and suppose that $A=\phi\,dS$ for
some smooth function $\phi$. Then using the standard fact that $\frac
{\partial}{\partial t}dS=-R\,dS$, we get%
\[
\frac{\partial}{\partial t}A=\left(  \frac{\partial}{\partial t}\phi
-R\phi\right)  dS.
\]
Because $\Delta_{d}A=\left(  \Delta\phi\right)  \,dS$, it follows that $A$
satisfies the Hodge heat equation $\frac{\partial}{\partial t}A=\Delta_{d}A$
if and only if $\phi$ evolve by equation (\ref{surface-system-a}). Now suppose
$E=-2df$ for some smooth function $f$. Then since%
\[
\left|  A\right|  ^{2}=\phi^{2}g^{ik}g^{j\ell}J_{i}^{p}g_{pj}J_{k}^{q}%
g_{q\ell}=2\phi^{2},
\]
it follows that $E$ satisfies $\frac{\partial}{\partial t}E=\Delta
_{d}E-d\left|  A\right|  ^{2}$ if and only if%
\[
-2d\left(  \frac{\partial}{\partial t}f\right)  =2\left(  d\delta+\delta
d\right)  \left(  df\right)  -4\phi\,d\phi=-2d\left(  \Delta f+\phi
^{2}\right)  .
\]
Hence we can apply Theorem \ref{main} with any solution $f$ of
(\ref{surface-system-b}).

To apply Corollary \ref{trace-Harnack}, we first calculate%
\[
\left(  \delta A\right)  _{i}=\left(  \delta\left(  \phi\,dS\right)
_{i}\right)  =\left(  \nabla^{j}\phi\right)  \left(  dS\right)  _{ij}%
=\nabla^{j}\phi J_{i}^{k}g_{kj}=J_{i}^{k}\nabla_{k}\phi=\left(  J\,d\phi
\right)  _{i}.
\]
Then we need only note that $\delta E=2\Delta f$ and set $V=2JX$ in the
resulting expression:%
\[
\psi\left(  A,E,V\right)  =\frac{1}{2}R\left|  V\right|  ^{2}-2\left\langle
\nabla\phi,JV\right\rangle +2\phi^{2}+2\Delta f.
\]

\end{proof}

\begin{proof}
[Proof of Example \ref{surface-trace}]The choice $\phi=tR+1$\ satisfies%
\[
\frac{\partial}{\partial t}\phi=t\left(  \Delta R+R^{2}\right)  +R=\Delta
\phi+R\phi.
\]
Then if $f=t\phi=t^{2}R+t$, we get%
\[
\frac{\partial}{\partial t}f=t\left(  \Delta\phi+R\phi\right)  +\phi=\Delta
f+\left(  tR+1\right)  \phi=\Delta f+\phi^{2}.
\]
Hence the trace LYH\ quadratic takes the form%
\[
\psi\left(  A,E,tX\right)  =t^{2}\left(  R\left\vert X\right\vert
^{2}+2\left\langle \nabla R,X\right\rangle +\Delta R+R^{2}\right)  +2tR+1,
\]
which is certainly positive at $t=0$.
\end{proof}

\begin{remark}
\label{strictly-positive}On any surface $\left(  \mathcal{M}^{2},g\left(
t\right)  \right)  $ of strictly positive curvature evolving by the Ricci
flow, the trace inequality
\[
0\leq R\left|  X\right|  ^{2}+2\left\langle \nabla\phi,X\right\rangle
+\frac{\partial}{\partial t}f
\]
in Proposition \ref{surface-system} can be deduced from the following
calculations. Let $F\doteqdot\Delta f+\phi^{2}-R^{-1}\left|  \nabla
\phi\right|  ^{2}$ and observe that%
\begin{align*}
\frac{\partial}{\partial t}F  &  =\Delta\left(  \Delta f+\phi^{2}\right)
+RF+R\phi^{2}+2\phi\Delta\phi-\left|  \nabla\phi\right|  ^{2}\\
&  -2R^{-1}\left(  \left\langle \nabla\phi,\nabla\Delta\phi\right\rangle
+\phi\left\langle \nabla\phi,\nabla R\right\rangle \right)  +R^{-2}\left|
\nabla\phi\right|  ^{2}\Delta R.
\end{align*}
Since by Bochner,%
\begin{align*}
\Delta F  &  =\Delta\left(  \Delta f+\phi^{2}\right)  -2R^{-1}\left(
\left\langle \nabla\Delta\phi,\nabla\phi\right\rangle +\left|  \nabla
\nabla\phi\right|  ^{2}+\frac{R}{2}\left|  \nabla\phi\right|  ^{2}\right) \\
&  +R^{-2}\left(  \left|  \nabla\phi\right|  ^{2}\Delta R+2\left\langle \nabla
R,\nabla\left|  \nabla\phi\right|  ^{2}\right\rangle \right)  -2R^{-3}\left|
\nabla R\right|  ^{2}\left|  \nabla\phi\right|  ^{2}.
\end{align*}
we have%
\begin{align*}
\frac{\partial}{\partial t}F  &  =\Delta F+RF+R\phi^{2}+2\phi\Delta\phi\\
&  +2R^{-1}\left(  \left|  \nabla\nabla\phi-R^{-1}\nabla R\otimes\nabla
\phi\right|  ^{2}-\phi\left\langle \nabla\phi,\nabla R\right\rangle \right) \\
&  \geq\Delta F+RF+R^{-1}\left(  \Delta\phi-\left\langle \nabla\phi,\nabla\ln
R\right\rangle +\phi R\right)  ^{2}.
\end{align*}
The inequality on the last line is equivalent to the estimate%
\begin{align*}
N  &  \doteqdot R\left(  R\phi^{2}+2\phi\Delta\phi\right)  +2\left(  \left|
\nabla\nabla\phi-R^{-1}\nabla R\otimes\nabla\phi\right|  ^{2}-\phi\left\langle
\nabla\phi,\nabla R\right\rangle \right) \\
&  -\left(  \Delta\phi-\left\langle \nabla\phi,\nabla\ln R\right\rangle +\phi
R\right)  ^{2}\\
&  =2\left(  \left|  \nabla\nabla\phi\right|  ^{2}-\left\langle \nabla\left|
\nabla\phi\right|  ^{2},\nabla\ln R\right\rangle +\left|  \nabla\phi\right|
^{2}\left|  \nabla\ln R\right|  ^{2}\right) \\
&  -\left(  \left(  \Delta\phi\right)  ^{2}-2\Delta\phi\left\langle \nabla
\phi,\nabla\ln R\right\rangle +\left\langle \nabla\phi,\nabla\ln
R\right\rangle ^{2}\right) \\
&  =2\left|  \nabla\nabla\phi-\nabla\phi\otimes\nabla\ln R\right|
^{2}-\left(  \Delta\phi-\left\langle \nabla\phi,\nabla\ln R\right\rangle
\right)  ^{2}\geq0.
\end{align*}

\end{remark}

\bigskip

Our final example makes no use of a K\"{a}hler structure:

\begin{example}
If $\left(  \mathcal{M}^{n},g\left(  t\right)  \right)  $ is any solution of
the Ricci flow on a closed manifold, the choice $A\equiv0$ lets us take
$E=-df$ for any solution $f$ of the ordinary heat equation $\frac{\partial
}{\partial t}f=\Delta f$, and leads to the LYH quadratic
\[
\Psi\left(  0,-df,U,W\right)  =\operatorname*{Rm}\left(  U,U\right)  +\left(
\nabla\nabla f\right)  \left(  W,W\right)  .
\]
In case the curvature operator of $\left(  \mathcal{M},g\right)  $ is
non-negative initially (hence for all time), Theorem \ref{main} applies
whenever $\nabla\nabla f\geq0$ initially and is equivalent to the statement
that the weak convexity of $f$ is preserved. (Compare Remark 2 in \S 6 of
\cite{CH}.)
\end{example}

\end{document}